\documentclass[12pt]{article}
\usepackage{amsmath}
\usepackage{amssymb}
\usepackage{amsfonts}
\newtheorem{thm}{Theorem}[section]
\newtheorem{pro}[thm]{Proposition}
\newtheorem{cor}[thm]{Corollary}
\newtheorem{lem}[thm]{Lemma}
\newtheorem{dfn}[thm]{Definition}
\newtheorem{fct}[thm]{Fact}
\newtheorem{note}[thm]{Note}
\newtheorem{rmk}[thm]{Remark}

=22cm\headheight=0cm\topskip=0cm\topmargin=0cm
\def\a{\bar {a}}\def\b{\bar {b}}
\def\x{\bar {x}}\def\y{\bar {y}}

\def\Ac{\mathcal {A}}
\def\Bc{\mathcal {B}}

\def\Rn{\mathbb R}

\def\qed{\hfill$\Box$}

\newcommand{\M}{\sf M}
\newcommand{\N}{\sf N}

\begin{document}
\def\dis{\displaystyle}

\begin{center}
{\Large{ Amenability, extreme  amenability,
   model-theoretic stability, \\ and dependence property in integral logic}} \vspace{10mm}

{\large{\bf Karim Khanaki}} \vspace{3mm}

{\footnotesize
  Faculty of Fundamental Sciences, Arak University of Technology,
 \\
P.O. Box 38135-1177, Arak, Iran; e-mail: khanaki@arakut.ac.ir
\\ \bigskip School of Mathematics,
Institute for Research in Fundamental Sciences (IPM), \\ P.O. Box
19395-5746, Tehran, Iran;
 e-mail: khanaki@ipm.ir} \vspace{5mm}
\end{center}

\begin{abstract}
This paper has three parts. First, we study and characterize
amenable  and extremely amenable topological semigroups in terms
of invariant measures using integral logic. We prove definability
of some properties of a topological semigroup such as amenability
and the fixed point on compacta property. Second, we define
types  and develop local stability in the framework of integral
logic. For a stable formula $\phi$, we prove definability of all
complete $\phi$-types over models and deduce from this the
fundamental theorem of stability. Third, we study an important
property in measure theory, Talagrand's stability. We point out
the connection between Talagrand's stability and dependence
property (NIP), and prove a measure theoretic version of
definability of types for NIP formulas.
\end{abstract}



{\small \tableofcontents

\bigskip
{\bf References} {\hspace{\stretch{0.01}} \bf \pageref{ref}} }

\section{Introduction} \label{}
Probability logics are logics of probabilistic reasoning. A model
theoretic approach aiming to study probability structures by
logical tools was started by   Keisler and   Hoover (see
\cite{Hoover,Keisler} for a survey). Among several variants of
this logic, they introduced integral logic $L_{\int}$ as an
equivalent `Daniell integral' presentation for $L_{\omega_1P}$.
Integral  logic uses the language of measure theory, i.e.,  that
of measurable functions and integration. The resulting framework
is close to the usual language of probability theory and allows
the formalization of much of probability. In \cite{BP} Bagheri
and Pourmahdian developed a finitary version of integration logic
and proved appropriate versions of the compactness theorem and
elementary JEP/AP. The intended models are graded probability
structures introduced by Hoover in \cite{Hoover} and in addition
to random variables over probability spaces, they include
dynamical systems and other interesting structures from real
analysis. In \cite{KB} the authors showed that many interesting
notions such as probability independence,  martingale property,
and  some special cases the notion of conditional expectation (as
in martingales) are expressible. Also, the Kolmogorov's extension
theorem was deduced from the compactness property of model
theory.  In \cite{KA} the authors further used the logical tools
to study invariant measures on compact Hausdorff spaces.
Consequently, they gave two proofs  of the existence of Haar
measure on compact groups. One might therefore hope to obtain
other applications of the compactness theorem.

Historically one of the great successes of model theory has been
Shelah's stability theory. Essentially the success of the program
largely due to the fact that certain (local) combinatorial
properties of formulas determine the corresponding global
properties. On the other hand, a general trend in  model theory
is to generalize these model-theoretic notions and tools to
frameworks that go beyond that of first order logic and
elementary classes.

In the present paper, on one hand, we study some analytic
concepts, amenability and extremely amenability, using integral
logic. On the other hand, we study  types and  local stability in
this logic. This approach has two advantages. First, we underline
the strengths of application of logical methods to the other
fields of mathematics.  Second, the results obtained by these
methods provide a new view on the related subjects in Analysis
and Logic, and open some fruitful areas of research on the
similar questions.

To summarize the results of this paper, in the first part
(Section \ref{4}), we consider an arbitrary topological semigroup
$S$ and any compact Hausdorff space $X$ such that $S$ acts
continuously on $X$ from the left. Let ${\rm{Inv}}_X(S)$ be the
set  of all Radon probability measures on $X$ which are left
invariant under elements of $S$. It is shown that the
nonemptiness of ${\rm{Inv}}_X(S)$ is expressible by a theory
$T_{S,X}$ in integral logic. We then present a characterization
of amenable topological semigroups in terms of invariant measures
(Fact~\ref{amenable-equiv}). Using the compactness theorem, we
give a proof of the fundamental result that goes back to
N.N.~Bogolioubov and N.M.~Krylov (Theorem~\ref{Krylov}).
 The interesting fact is that for a topological semigroup $S$ the amenability of
  $S$ is expressible by a theory $T_S$ in the framework of integral logic. Some other
new results and different proofs of some known results are given
for extremely amenable topological semigroups (Fact~\ref{2}, and
Propositions~\ref{extr-amen}).

Although most of the results in the first part of the paper are
standard, the study of amenable and extremely amenable semigroups
is necessary because it leads us to the ``true and correct" notion
of a type  in integral logic. In fact, types are known
mathematical objects, Riesz homomorphisms. Thus, for a complete
theory $T$, the space of complete types $S(T)$ can be represented
by the {\em spectrum} of $T$.  Thereby, in the second part of the
paper (Section \ref{6}), we define types and develop local
stability. For a stable formula $\phi$, we prove that all
complete $\phi$-types over models are definable, and we deduce
from this the fundamental theorem of stability
(Corollary~\ref{Fundamental-Thm}). We show that a formula $\phi$
is stable if and only if  its Cantor-Bendixson rank is finite.

 In the third part of the paper (Section \ref{NIP}),
  we study a form of the dependence property which is
an important measure-theoretic property, {\em Talagrand's
stability}. Then we prove that for an {\em almost dependent
formula} $\phi$, all $\phi$-types are {\em almost} definable
(Theorem~\ref{almost-Definable-Thm}). We then study the
Cantor-Bendixson rank in almost dependent theories.

It is  worth recalling another line of research arisen from ideas
of Chang and Keisler \cite{CK}, namely continuous logic. The idea
was recently refined and developed in \cite{BU} and \cite{BBHU} by
Ben~Yaacov, Berenstein, Henson, and  Usvyatsov for the class of
metric structures which include such important classes of
structures as Banach spaces and measure algebras. Although some
results in the present paper (cf. Section~\ref{6}) are similar to
those in \cite{BU}, in some senses they are different: (i) Our
approach can be used to generalize the results in \cite{BU} and
\cite{Mofidi} (see Remark \ref{topometric}); (ii) In \cite{Ben1}
and \cite{Ben2}, Ben~Yaacov proved that the theory $ARV$ and the
category of probability algebras are $\aleph_0$-stable. Note that
in this paper we do not study probability measure algebras or
$L^1$-spaces, but we study measurable functions.  In contrast to
\cite{Ben1} and \cite{Ben2}, the theory of a probability
structure is not necessarily stable. This leads us to the
dichotomy between stable probability structures and unstable
probability structures; (iii) Some analytic properties such as
probability independence, amenability, extreme amenability and
the existence of invariant measures on compact spaces are
expressible in the framework of integral logic.

After the submission of the present paper we came to know that,
independently from us, similar ideas were used by P. Simon
\cite{S} in classical logic. We note that the argument for almost
definability in the case of a dependent formula is truly measure
theoretic and can be used for proving some new results in
classical logic. We will study it in a future work.

The organization of the paper is as follows. In the next section
we review some basic notions from measure theory. In
Section~\ref{integral-logic}  a summary of results on integral
logic from \cite{BP} are given. In Section \ref{4}, we study
amenable and extremely amenable topological semigroups, and give
a characterization of (extreme) amenability in terms of
(multiplicative) invariant measures. A proof of the
Bogolioubov-Krylov theorem is given in Section~\ref{4}.  It is
shown that the (extreme) amenability of a topological semigroup
$S$ is expressible by a theory $T_S$ ($\texttt{T}_S$) within
integral logic. In Section~\ref{6}, we conclude with the
development of local stability, and we prove the fundamental
theory of stability. In Section~\ref{NIP}, we study NIP theories
and give some results.

\section{Preliminaries from topological measures theory} \label{sec2}
In this section we review some basic notions from measure theory.
Further details can be found in \cite{Folland,Fremlin2,Fremlin4}.
 Let $X$ be a compact Hausdorff space.
 The space $C(X,{\Bbb R})$ of continuous
real-valued functions on $X$  is denoted here by $C(X)$.  Since
$X$ is a compact space, every $f\in C(X)$ is bounded and $C(X)$
is a normed vector space with the uniform norm.

 The class of {\em Baire sets} is defined to be the
smallest $\sigma$-algebra $\cal B$ of subsets of $X$ such that
each function in $C(X)$ is measurable with respect to $\cal B$.
 The smallest $\sigma$-algebra containing the open sets is called
the class of {\em Borel sets}. Clearly, every
 Baire set is a Borel set, but there are
compact spaces where the class of Borel sets is larger than the
class of Baire sets. 
 By a {\em Baire} ({\em Borel}) {\em measure} on $X$ we mean a finite measure
 defined for all Baire (Borel) sets.
A {\em Radon measure} on $X$ is a Borel measure which is regular.
It is known that every Baire measure on a compact space is regular
and has a unique extension to a Radon
 measure. 

 A {\em topological semigroup} is a semigroup $S$ endowed
with a Hausdorff topology such that the operation $(x, y) \mapsto
xy$ is continuous from $S \times S$ to $S$.
 By a {\em topological group} we mean a group $G$
endowed with a Hausdorff topology such that the group operations
$(x, y) \mapsto xy$ and $x \mapsto x^{-1}$ are continuous from
$G\times G$ and $G$ to $G$. A  topological group whose topology
is (locally) compact and Hausdorff is called a ({\em locally})
{\em compact} {\em group}.

A  topological semigroup $S$ is said to {\em act on a topological
space $X$ from the left} if there is a map $S\times X\to X$
(denoted by $(s, x)\mapsto s\cdot x$   for each $(s, x) \in
S\times X$) such that (a) the map $x\mapsto s\cdot x$ is
continuous for each $s \in S$, (b) for $s, s'\in S$, $(ss')\cdot
x = s\cdot (s'\cdot x)$ for each $x\in X$, and (c) if $S$ has the
identity $e$, then $e\cdot x=x$ for each $x\in X$. In addition,
the left action is said to be {\em continuous} if $(s, x)\mapsto
s\cdot x$ is a continuous map from $S\times X$ to $X$. Similarly
one can define a right (continuous) action.  If $S$ acts on
topological space $X$ from the left (right) and $E\subseteq X$
and $s\in S$ we define
$$s\cdot E=\{s\cdot x:x\in E\} \ \ \ \ \ \ (E\cdot s=\{x\cdot s:x\in E\}).$$

If $f$ is a continuous real-valued function on a topological space
$X$ and $s\in S$, we define the {\em left} ({\em right}) {\em
translate of $f$ by} $s$, as follows:
$$(f\cdot s)(x)=f(s\cdot x) \ \ \ \ \ \ ((s\cdot f)(x)=f(x\cdot s)).$$
The point of the above definition is to make $f\cdot
(ss')=(f\cdot s)\cdot s'$~~($(ss')\cdot f=s\cdot(s'\cdot f)$).

If a topological semigroup $S$ acts on a space $X$ from the left
(right), a measure $\mu$ on $X$ is {\em left} ({\em right}) {\em
$S$-invariant}  if $\mu(s\cdot E)$  ($\mu(E\cdot s)$) is defined
and equal to $\mu(E)$ whenever $s\in S$ and $\mu$ measures $E$.
If $X$ be   a compact Hausdorff space, then a linear functional
$I$ on $C(X)$ is called {\em left} ({\em right}) {\em
$S$-invariant} if $I(f\cdot s)=I(f$)  ($I(s\cdot f)=I(f)$) for
all $s$ in $S$ and $f$ in $C(X)$.

A {\em left} ({\em right}) {\em Haar measure} on a compact group
$G$ is a nonzero left (right) $G$-invariant  Radon measure $\mu$
on $G$.

\begin{pro}  \label{1}  {\em (\cite[Proposition~441L]{Fremlin4})}
Let $X$ be a  Hausdorff  compact space and $S$ a topological
semigroup which acts on $X$. A nonzero Radon measure $\mu$ on $X$
is a left (right) $S$-invariant measure iff $\int fd\mu=\int
(f\cdot s) d\mu$ ($\int fd\mu=\int (s\cdot f)d\mu$) for all $f\in
C(X)$ and $s\in S$.
\end{pro}

If $G$ is compact group, then a left Haar measure on $G$ is also
a right Haar measure. Also, the Haar measure is unique up to a
positive scalar multiple, i.e.  if $\mu$ and $\nu$ are Haar
measures on a compact group $G$, there exists $c>0$ such that
$\mu=c\nu$.

\bigskip
\noindent
\textbf{The Riesz Representation Theorem.} \label{Riesz}
 {\em Let $X$ be a locally compact Hausdorff space and $C_c(X)$ the space of continuous real-valued functions on $X$
with compact support.
\begin{itemize}
\item [{\em (a)}] {\em (\cite[p. 212]{Folland})} If $I$ is a positive linear functional on
  $C_c(X)$, there is a unique Radon measure $\mu$ on
$X$ such that $I(f) = \int f d\mu$ for all $f\in C_c(X)$.
\item [{\em (b)}] {\em (\cite[p. 358]{Royden})} If $X$ is compact, then the dual
of $C(X)$ is (isometrically isomorphic to) the space of all finite
signed Baire measures on X with norm defined by
$\|\mu\|=|\mu|(X)$.
\end{itemize} }

\noindent \textbf{The Hahn-Banach Theorem.} \label{Hahn-Banach}
{\em {\em (\cite[p. 159]{Folland})} Let $\cal N$ be a normed
vector space. If $\cal M$ is a closed subspace of $\cal N$ and $x
\in {\cal N}\setminus{\cal M}$, there exists a bounded linear
functional $I$ on $\mathcal{N}$ such that $I|_{\cal M}=0$,
$\|I\|=1$ and $I(x)=\inf_{y\in{\cal M}}\|x-y\|$. }

\bigskip Let $(M,\Bc,\mu)$ be a measure space and $\mu^*$ its
associated outer measure defined by
$$\mu^*(X)=\inf\{\mu(A)|\ X\subseteq A\in\Bc\}.$$
If $N\subseteq M$, then $\Bc_N=\{A\cap N|\ A\in \Bc\}$ is a
$\sigma$-algebra and $\mu_N=\mu^*\upharpoonright\Bc_N$ is a
measure on $N$. $\mu_N$ is called the {\em subspace measure} on
$N$. A measurable envelope for $N$ is a measurable set $E\in\Bc$
such that $N\subseteq E$ and $\mu(E\cap A)=\mu^*(N\cap A)$ for any
$A\in\Bc$. Every $N\subseteq M$ of finite outer measure has an
envelope (e.g take $E\in\Bc$ containing $N$ with
$\mu(E)=\mu^*(N)$). If $f:M\rightarrow \Rn$ is measurable,
$\int_N f$ abbreviates $\int_N (f\upharpoonright N)d\mu_N$.

\begin{pro} \label{submeasure} {\em (\cite[p. 38]{Fremlin2})}
Let $(M,\Bc,\mu)$ be a measure space, $N\subseteq M$, and $f$ be
an integrable function defined on $M$.
\begin{itemize}
\item [{\em (a)}] If $f$ is nonnegative then $f\upharpoonright N$ is $\mu_N$-integrable and
$\int_N f\le\int f$.
\item [{\em (b)}] If either $N$ is of full outer measure in $M$ or
$f$ is zero almost everywhere on $M-N$, then $\int_N f=\int_M f$.
\end{itemize}
\end{pro}

\section{Integral logic} \label{integral-logic}
In this section we give a brief review of integral logic from
\cite{BP,KA}. Results from \cite{BP,KA} are stated without proof.
All languages are assumed to contain  unary relation and constant
symbols.  Let $L$ be a language. To each relation symbol $R\in L$
we assign a nonnegative real number
 $\flat_R\geqslant 0$ called
the universal bound of $R$. The terms are just the constant
symbols and the variables.

\begin{dfn} \label{formula}
\em{The family of $L$-formulas and their universal bounds is
defined as follows:
\begin{itemize}
  \item [1.] If $R$ is a relation  symbol and $t$ is a term, then
  $R(t)$ is an atomic formula with bound $\flat_R$.
  \item [2.] If $\phi$ and $\psi$ are formulas and
  $r,s\in\mathbb{R}$, then so are $r\phi+s\psi$ and $\phi\times\psi$
  with bounds  $|r|\flat_\phi+|s|\flat_\psi$ and
  $\flat_\phi\flat_\psi$, respectively.
  \item [3.] If $\phi$ is a formula, then $|\phi|$
  is a formula with bound $\flat_\phi$.
  \item [4.] If $\phi$ is a formula and $x$ is a variable, then $\int\phi dx$
  is a formula with bound $\flat_\phi$.
\end{itemize}}
\end{dfn}

Note that $\phi^+=\frac{1}{2}(\phi+|\phi|)$ and
$\max(\phi,\psi)=(\phi-\psi)^++\psi$ and  similarly $\phi^-$ and
$\min(\phi,\psi)$ are formulas.

\begin{dfn} \label{structure}
{\em An {\em  $L$-structure} is a probability measure space
${\M}=(M,{\cal B},\mu)$ equipped with:
\begin{itemize}
  \item  for each constant symbol $c\in L$, an element $c^{\M}\in M$;
  \item  for each relation symbol $R\in L$, a measurable
 map $R^{\M}:M\rightarrow[-\flat_R,\flat_R]$.
\end{itemize}}
\end{dfn}

$L$-structures are denoted by $\M,\N$ etc. The notion of free
variable is defined as usual and one writes $\phi(\x)$ (or
$\phi(x_1,\ldots,x_n)$) to display them. If $\M$ is an
$L$-structure, for each formula $\phi(x_1,\ldots,x_n)$ and $\a\in
M^n$, $\phi^{\M}(\a)$ is defined inductively starting from atomic
formulas. In particular,
$$\bigg(\int\phi(\x,y)dy\bigg)^{\M}(\a)=\int\phi^{\M}(\a,y)dy.$$ An
easy induction  shows that every $\phi^{\M}(\x)$ is a
well-defined measurable function from $M^n$ to
$[-\flat_{\phi},\flat_{\phi}]$. Indeed, for every $\phi(\x,\y)$
and $\a$, $\phi^{\M}(\a,\y)$ is measurable. Moreover, we have
$\int\int\phi dxdy=\int\int\phi dydx$.

A formula is closed if no free variable occurs in it. A
\emph{statement} is an expression of the form $\phi(\bar{x})\geq
r$ or $\phi(\bar{x})=r$. Closed statements are defined similarly.
Any set of closed statements is called a theory. The theory of a
structure $\sf M$ is the collection of closed statements
satisfied in it. Such theories are called complete. $\M,\N$ are
{\em elementarily equivalent} (written $\M\equiv \N$) if they
have the same theory. The notion ${\M}\vDash \Gamma$ is defined
in the obvious way. If $T$ is an $L$-theory, two formulas
$\phi(\bar{x}),\psi(\bar{x})$ are said to be $T$-equivalent if
the statement $\phi=\psi$   a.e. is satisfied in every model of
$T$. We say $T$ has quantifier-elimination if every formula is
$T$-equivalent to a quantifier-free formula (i.e. without $\int$).

The ultaproduct of a family ${\M}_i$, $i\in I$ of structures over
an ultrafilter $\mathcal D$ is an $L$-structure and denoted by
${\M}=\prod_{\mathcal D}{\M}_i$ (cf. \cite{BP,KA}).

\begin{thm} {\em (Fundamental theorem)} \label{fundamental}
For each $\phi(\x)$ and $[a_i^1],\ldots,[a_i^n]\in M$
$$\phi^{\M}([a_i],\ldots,[a_i^n])=\lim_{\mathcal D}\phi^{{\M}_i}(a^1_i,\ldots,a_i^n).$$
\end{thm}

\bigskip
An immediate consequence of the fundamental theorem is the
following whose proof is just a modification of its analog in the
usual first order logic.

\begin{thm} {\em (Compactness theorem)} \label{compact-thm}
Any finitely satisfiable set of closed statements is satisfiable.
\end{thm}

\begin{dfn}
{\em (i) If $M\subseteq N$, $\M$ is a \emph{substructure} of $\N$,
denoted by ${\M}\subseteq{\N}$, if $\M$ has the subspace measure
and for each $R\in L$ and $\a\in M$, $R^{\M}(\a)=R^{\N}(\a)$. If
these equalities hold for almost all $\a$, $\M$ is called an
\emph{almost substructure} of $\N$ and is denoted by
${\M}\subseteq_a{\N}$.

(ii) An injection $f:\M\rightarrow \N$ is called an
\emph{elementary embedding} if for each $\phi$ and $\a\in M$,
$\phi^{\M}(\a)=\phi^{\N}(f(\a))$. It is an \emph{almost
elementary embedding} if for each $\phi$ this holds almost surely
for $\a\in M$.  If $f$ is the inclusion, these are respectively
denoted by ${\M}\preceq {\N}$, and ${\M}\preceq_a {\N}$. $f$ is
said to be {\em almost surjective} if its range has full measure.
One also defines \emph{isomorphism} (resp. \emph{almost
isomorphism}) as a surjective  (resp. almost surjective)
elementary (resp. almost elementary) embedding.}
\end{dfn}

The fact that
 $\preceq$ (resp. $\preceq_a$) is stronger than
$\subseteq$ (resp. $\subseteq_a$) is a consequence of the
Tarski-Vaught test (see below). Among the two notions of
isomorphism, the notion of almost isomorphism is more useful for
us, however, the exact isomorphism appears naturally in some
cases. In ergodic theory, a map which is an (exact) isomorphism
after removing some negligible sets from its domain and codomain
is called an isomorphism. This notion is equivalent to our notion
of almost isomorphism.

A  structure is called \emph{minimal} if it has no redundant
measurable sets, i.e., for any substructure
${\M}'=(M,\Ac,\mu\upharpoonright \Ac)$ where $\Ac\subseteq\Bc$,
one has $\Ac=\Bc$. In fact, every structure is isomorphic to a
minimal structure, which can be explicitly described.

\begin{pro} \label{minimal}
Let ${\M} = (M,{\cal B},\mu)$ be an $L$-structure and $\cal A$ be
the $\sigma$-algebra generated by the sets of the form $\{x\in
M:\phi^{\M}(x)>0\}$ where $\phi$ is any formula with parameters in
$M$. Then, ${\M}' = (M,{\cal A},\mu\upharpoonright{\cal A})$ is a
minimal measure $L$-structure isomorphic to $\M$.
\end{pro}

\begin{pro}[Tarski-Vaught Test for $\preceq$] \label{Tarski-Vaught}
Let ${\M}, {\N}$ be minimal. If $M\subseteq N$ then $\M\preceq
\N$ if and only if for each $\phi(\a,x)$, where $\a\in M$, the
intersection of the set $\{x\in N: \ \phi^{\N}(\a,x)>0\}$ with
$M$ is $\mu_{\M}$-measurable and has the same measure. Similar
statement holds for ${\M}\preceq_a {\N}$ with `for almost all
$\a$' in place of `for each $\a$'. In both cases,
$\mu_{\M}=\mu_{\N}\upharpoonright M$.
\end{pro}

Next we are going to prove a key result, which plays an important
role in the rest of this paper. Assume that $X$ is a compact
Hausdorff space.  Let $L_X$ be the language
 consisting of a unary relation symbol $R_f$ for each $f\in
 C(X)$ and a constant symbol $c_a$ for each $a\in X$.  Let ${\M}$ be an
 $L_X$-structure with the following properties:
\begin{itemize}
  \item   $X\subseteq{M}$;
  \item   the restriction of  $R_f^{{\M}}$ to $X$ is $f$, particularly $R_1^{\M}={\bf 1}$;
  \item   $R_{f+g}^{\M}=R_f^{\M}+R_g^{\M}$ and  $R_{r\times f}^{\M}=r\cdot R_f^{\M}$ for each $f,g\in C(X)$ and real number $r$;
  \item   $R_{f\times g}^{\M}=R_f^{\M}\times R_g^{\M}$ for each $f,g\in C(X)$;
  \item   $R_{\max(f,g)}^{\M}=\max(R_f^{\M},R_g^{\M})$ for each $f,g\in C(X)$.
\end{itemize}

The next proposition shows that the subspace  measure $\mu_X$ on
$X$ behaves like the measure $\mu$ on ${\M}$. In fact, $(X,{\cal
B}_X,\mu_X)$ with the natural
 interpretation of relation and constant symbols
 is an elementary substructure of $\M$.

\begin{pro} \label{Pushing-Down}
Assume that $X$ and ${\M}$ are as above.
\begin{itemize}
\item [{\em (a)}] The subspace measure $\mu_X$ on $X$ is a regular Baire
measure such that $\int fd\mu_X=\int R_f^{\M}d\mu$ for each $f\in
C(X)$.
\item [{\em (b)}] There exists a  Radon  measure $\bar{\mu}_X$ on
$X$ such that $\int fd\bar\mu_X=\int R_f^{\M}d\mu$ for each $f\in
C(X)$.
\end{itemize}
\end{pro}
{\bf Proof.} (a).  By Proposition \ref{submeasure}, it suffices to
show that $X$ is of full outer measure in $M$. We assume that
${\M}$ is minimal. By Proposition \ref{minimal},
$$\mu_X(X)=\inf\Bigg\{\sum_1^\infty\mu(A_k):X\subseteq\bigcup_1^{\infty} A_k\Bigg\}$$
where $A_k=(R_{f_k}^{\M})^{-1}(0,\infty)$ for a  $f_k\in C(X)$
because every formula $\phi$ is equal to a relation symbol $R_f$.
We show that $\mu(\bigcup_k A_k)=1$ for every sequence $\langle
A_k\rangle_{k\in \mathbb{N}}$ such that
$X\subseteq\bigcup_1^{\infty} A_k$.  If $X\subseteq \bigcup_k
f_k^{-1}(0,\infty)$, then there exist $f_1,\dots,f_n$ such that
$X=\bigcup_1^n f_k^{-1}(0,\infty)$ because $X$ is compact. If
$f=\max(f_1,\ldots,f_n)$, then $X=f^{-1}(0,\infty)$. Thus,
$X\subseteq (R_f^{\M})^{-1}(0,\infty)$ because
$R_f^{\M}=\max(R_{f_1}^{\M},\ldots,R_{f_n}^{\M})$. Since $X$ is
compact and $f$ is continuous, there exist real numbers
$s\geqslant r>0$ such that $X=f^{-1}[r,s]$. Also, we can easily
check that $M=(R_f^{\M})^{-1}(0,\infty)$ since $R_f^{\M}\geqslant
r$.  Thus, $\mu(\bigcup_k
A_k)\geqslant\mu((R_f^{\M})^{-1}(0,\infty))=1$, i.e.
$\mu_X(X)=1$. We may assume that $\mu_X$ is a Baire measure.
Also, we know that every Baire measure on a compact space is
regular.

(b).  It is known that every Baire regular measure on a compact
space has a unique extension to a Radon measure (cf. \cite[p. 341]
{Royden}). Let $\bar\mu_X$ be the unique extension of $\mu_X$ to
a Radon measure on $X$. Since only the values of $\bar\mu_X$ on
 Baire sets matter for $\int fd\bar\mu_X$, we
 have
$\int fd\bar\mu_X=\int fd\mu_X$ for each $f\in C(X)$.~\qed

\section{Amenability and extreme amenability} \label{4}
In this section  we study and characterize amenable and extremely
amenable  topological semigroups in terms of invariant measures
using integral logic. First, we give two conditions equivalent to
the existence of measures on a compact Hausdorff space $X$
invariant under a semigroup $S$ which acts on it from the left.
We then characterize (extremely) amenable topological semigroups
in terms of (multiplicative) invariant measures.  It is shown that
all compact groups, abelian topological semigroups, and all
locally finite topological groups are amenable. An interesting
fact is that for a topological semigroup $S$ the (extreme)
amenability of $S$ is expressible by a theory $T_{S}$
(${\texttt{T}}_S$) in the framework of integral logic. Therefore,
it is shown that a locally compact group $G$ has no Borel
paradoxical decomposition iff the theory $T_{G}$ is satisfiable.

  Let $X$ be a compact  Hausdorff space and $S$ be a semigroup which acts on $X$ from the left.
Let $L_X$ be the language consisting of a unary relation symbol
$R_f$ for each $f\in
 C(X)$ and a constant symbol $c_a$ for each $a\in X$
and $T_{S,X}$ be the theory with the  following  axioms:
\begin{itemize}
\item [(1)] $R_{1}={\bf  1}$,
\item [(2)] $\int R_1dx=1$,
\item [(3)] $ R_{f}(c_a)=f(a)$ \ \ \ \ \ \ \ \ \ \ \ \ \ \ \ \ \ \ \ \ \ \ \ \ \ \ \ \ \ for each $R_f,c_a\in L_X$,
\item [(4)] $ R_{f+g}=R_f+R_g$ \ \ \ \ \ \ \ \ \ \ \ \ \ \ \ \ \ \ \ \ \ \ \ \ ~ for each $R_f,R_g\in L_X$,
\item [(5)] $R_{r\times f}=r\times  R_f$ \ \ \ \ \ \ \ \ \ \  \ \ \ \ \ \ \ \ \ \ \ \ \ \ \ \ \ \  for each $R_f\in L_X$ and   $r\in \mathbb{R}$,
\item [(6)] $ R_{f\times g}=R_f\times R_g$ \ \ \ \ \ \ \ \ \ \ \ \ \ \ \ \ \ \ \ \ \ \ \ \ \ ~ for each $R_f,R_g\in L_X$,
\item [(7)] $ R_{\max(f,g)}=\max(R_f,R_g)$ \ \ \ \ \ \ \ \ \ \ \ \ \ \ \ ~ for each $R_f,R_g\in L_X$,
\item [(8)] $ \int R_f(x)dx=\int R_{(f\cdot s)}(x)dx$ \ \ \ \ \ \ \ \ \ \ \ \ \ ~  for each $R_f\in L_X$ and  $s\in S$,

 where  $(f\cdot s)(x)=f(s\cdot x)$.
\end{itemize}
 \vspace{.3 cm}
 Note that (1)  says that the interpretation of $R_1$
 is the constant function $\bf 1$, (2) means that we have a probability measure,
 (3) says that $f$ is a subset of the interpretation of $R_f$,
 (4)$-$(7) that the family of the interpretations of
relation symbols is a vector lattice,  and (8) means that the
measure is left $S$-invariant. $T_{S,X}$ is called the theory of
{\em left $S$-invariant measures on $X$}.

As a  consequence  of the compactness theorem we give conditions
equivalent to the existence of  a left $S$-invariant Radon
measure on  $X$. Later, we give results based on these
conditions. Let ${\rm{Inv}}_X(S)$ be the set of all regular Borel
probability measures on $X$ which are
 left $S$-invariant.

\begin{pro} \label{Prop-1}
Assume that $S, X$ and $T_{S,X}$ are as above. Then the following
are equivalent:
\begin{itemize}
           \item [{\em (i)}] ${\rm{Inv}}_X(S)\neq\emptyset$.
           \item [{\em(ii)}] $T_{S,X}$ is  satisfiable.
\end{itemize}
\end{pro}
{\bf Proof.} (i)$\Rightarrow$(ii) is obvious.
 For the converse, let ${\M}$ be a model of $T_{S,X}$. By Urysohn's lemma,  one can easily
verify that $X\subseteq M$.
 By Proposition~\ref{Pushing-Down}(b), there exists a Radon measure
 $\bar\mu_X$ on $X$ such that $\int f d\bar\mu_X=\int R_f^{\M}
d\mu$ for each $f\in C(X)$. Therefore, $\bar\mu_X$ is a nonzero
regular Borel left $S$-invariant measure on $X$.~\qed

\bigskip
The following classical result gives a
 condition equivalent to the existence of a
left $S$-invariant Radon measure on $X$ (see \cite{Hewitt},
Theorem~17.15).

\begin{fct} \label{Compact-Haar-3}
Let $S$ be a semigroup with identity. If $S$ acts from the left
on a  compact Hausdorff space $X$, then the following are
equivalent:
\begin{itemize}
           \item [{\em (i)}] ${\rm{Inv}}_X(S)\neq\emptyset$.
           \item [{\em (ii)}]  For every elements $s_1,\ldots,s_n$ of $S$
and elements $f_1,\ldots,f_n$ of $C(X)$ we have
$$\bigg\| {\bf 1}-\sum_{i=1}^n(f_i\cdot s_i-f_i)\bigg\|\geq 1.$$
\end{itemize}
\end{fct}
{\bf Proof.} (i)$\Rightarrow$(ii): Let $h= \sum_{i=1}^n(f_i\cdot
s_i-f_i)$. If $\sup_{x\in X}|{\bf 1}-h(x)|=1-\epsilon$ where
$\epsilon$ is a positive real number, then $\epsilon<h(x)<2$ for
all $x\in X$, thereby $\int hd\mu>\epsilon$ for every probability
measure $\mu$ on $X$, i.e., ${\rm{Inv}}_X(S)=\emptyset$. For the
converse, let $L_X$ be the language consisting of a unary
relation symbol $R_f$ for each $f\in
 C(X)$ and a constant symbol $c_a$ for each $a\in X$
and $T_{S,X}$ be the theory of left $S$-invariant measures on
$X$.  By Proposition~\ref{Prop-1}, it suffices to show that
$T_{S,X}$ is finitely satisfiable. Assume that $\Gamma$ is a
finite subset of $T_{S,X}$ such that for each $i\leq n$ and
$j\leq m$ the statement $\int R_{f_i}dx=\int R_{f_i\cdot s_j}dx$
is in $\Gamma$. Thus, $f_1,\ldots,f_n$ are in $C(X)$ and
$s_1,\ldots,s_m$ are in $S$. Let $M$ be the closure of the
subspace generated by $f_i-f_i\cdot s_j$ for each $i\leq n$ and
$j\leq m$. Since $S$ has an identity, clearly $\inf_{h\in
M}\|{\bf  1}-h\|=1$. Let $K$ be a subspace of $C(X)$ such that
$M+K=C(X)$ and $M\cap K=0$. By the Hahn-Banach theorem, define
$I$ to be $0$ on $M$ and a nonzero bounded  linear functional on
$K$ such that $I({\bf 1})=\|I\|=1$.  By the Riesz representation
theorem, there exists a signed Baire measure $\mu$ on $X$ such
that $\int (f_i-f_i\cdot s_j) d\mu=0$
 for each $i\leq n$ and $j\leq m$. Also, $\mu$ is a nonzero positive measure because
$\mu(X)=\int {\bf  1}d\mu=I({\bf  1})=\|I\|=|\mu|(X)$. Hence
$(X,\mu)$ with the natural interpretation of relation and
constant symbols is a model of $\Gamma$.~\qed

\subsection{Amenability}
In this subsection we define amenable topological semigroups
 and characterize  them in
terms of invariant measures. Also, we show that all compact groups
and locally finite topological groups are amenable.
 Let $S$ be a topological semigroup,
  and   $C_b(S)$  the Banach
space of all bounded real-valued continuous functions on $S$ with
the usual supremum norm. For $s\in S$ and $f \in C_b(S)$, let
$f\cdot s$ and $s\cdot f$ be the elements in $C_b(S)$ defined by
$(f \cdot s)(t) = f(st)$ and $(s \cdot f)(t) = f(ts)$, $t\in S$,
respectively. A subspace $E$ of $C_b(S)$ is {\em left} ({\em
right}) {\em invariant} if $f\cdot s\in E$  ($s\cdot f\in E$) for
all $s\in S, f\in E$. If $E$ is  both left and right
 invariant, then $E$ is called {\em invariant}.

Let $E$ be a left invariant closed subspace of $C_b(S)$ that
contains $\bf 1$, the constant $1$ function on $S$.
 A {\em mean} on  $E$  is a linear functional $I$ on $E$ such that
\begin{itemize}
           \item [(1)] $I({\bf  1}) = 1$,
           \item [(2)] $I(f)\geq 0$ if $f\geq 0$.
\end{itemize}
 A mean $I$ on a  left (right)  invariant closed subspace $E$   of $C_b(S)$ that
contains $\bf 1$ is said to be {\em left} ({\em right}) {\em
invariant} if $I(f \cdot s) = I(f)$  ($I(s\cdot f) = I(f)$) for
all $f \in E$ and $s\in S$.

We define the subspace $LUC(S)$ of all {\em left uniformly
continuous functions} in $C_b(S)$ which plays an important role in
the rest of this paper.  For a topological semigroup $S$ set
 $$LUC(S)=\{f\in C_b(S):~\mbox{the
map } s\mapsto  f\cdot s \textrm{ is (norm) continuous from $S$ to
$C_b(S)$} \}.$$ Similarly one can define the subspace $RUC(S)$ of
all {\em right uniformly continuous functions} in $C_b(S)$. It is
known that $LUC(S)$ and $RUC(S)$  are closed and invariant
subalgebras  of $C_b(S)$. They are also closed under the lattice
operations (cf. \cite[Lemmas~1.1 and~1.2]{Namioka}). Therefore,
$LUC(S)$ and $RUC(S)$  are $M$-spaces with the unit $\bf 1$.

 \begin{dfn}
 {\em A topological semigroup
$S$ is said to be {\em left} ({\em right}) {\em amenable} if
$LUC(S)$ ($RUC(S)$) admits a left (right) invariant mean. A
topological semigroup $S$ is called {\em amenable} if it is both
left and right amenable.}
\end{dfn}

We now characterize amenable topological semigroups in terms of
invariant measures, for which we need the following lemma.

\begin{lem} \label{lemma}
Let $S$ be a topological semigroup.
\begin{itemize}
           \item [{\em (i)}] If  $X$ is a  closed and invariant
           subset of $\{I\in LUC(S)^*:\|I\|=1\}$, then  the natural action of $S$ on  $X$  is continuous.
           \item [{\em (ii)}]  If $X$ is a  compact Hausdorff
space and $\cdot$ is a  continuous action of $S$ on $X$ (by the
left side), then, for each $f\in C(X)$,
 the map $s\mapsto f\cdot s$ from $S$ to $C(X)$ is (norm)
continuous.
\end{itemize}
\end{lem}
{\bf Proof.}  (i): Assume that $s,s'\in S$ and $I,I'\in X$. Then
for each $f\in LUC(S)$  we have
\begin{align*}
|(s'\cdot I')(f)-(s\cdot I)(f)| &=|I'(f\cdot s')-I(f\cdot s)|  \\
&\leq |I'(f\cdot s')-I'(f\cdot s)|+|I'(f\cdot
s)-I(f\cdot s)| \\
&= |I'(f\cdot s' -f\cdot s)|+|I'(f\cdot
s)-I(f\cdot s)| \\
&\leq \|I'\|\times \|f\cdot s' -f\cdot s\| +|I'(s\cdot f)-I(s\cdot
f)|\\&= \|f\cdot s' -f\cdot s\| +|I'(f\cdot s)-I(f\cdot s)|.
\end{align*}
Therefore the continuity of $(s,I)\mapsto I\cdot s$ follows from
the continuity $s\mapsto f\cdot s$.

(ii): Let $f\in C(X)$, $s_0\in S$ and $\epsilon>0$, and let $U$ be
the subset of $S\times X$  given by $U=\{(s,x):|f(s_0\cdot
x)-f(s\cdot x)|<\epsilon\}$. Then $U$ is open and $\{s_0\}\times
X\subseteq U$.  Hence there is a neighborhood $V$ of $s_0$ such
that $V\times X\subseteq U$, and it follows that $\|f\cdot
s_0-f\cdot s\|<\epsilon$ whenever $s\in V$.~\qed

\bigskip
We now give a classical result.

\begin{fct} \label{amenable-equiv}
Let $S$ be a topological semigroup with identity. Then the
following are equivalent:
\begin{itemize}
           \item [{\em (i)}] $S$ is left  amenable.
           \item [{\em (ii)}] Whenever  $X$ is a non-empty compact Hausdorff
space and $\cdot$ is a  continuous action of $S$ on $X$ (by the
left side), then  ${\rm{Inv}}_X(S)\neq\emptyset$.
\end{itemize}
\end{fct}
{\bf Proof.} (i)$\Rightarrow$(ii): By Fact~\ref{Compact-Haar-3} it
suffices to show that $\sup_{x\in X}|{\bf 1}-h(x)|\geq 1$, which
$h$ is of the form $\sum_{i=1}^n(f_i\cdot s_i-f_i)$ where
$s_1,\ldots,s_n$ are elements of $S$ and $f_1,\ldots,f_n$ are in
$C(X)$. If not, then $\sup_{x\in X}h(x)<0$. Let $I$ be  a left
invariant mean on $LUC(S)$.
 Fix a positive linear functional $\Lambda$ on $C(X)$.
Define $\widetilde{f}:S\to \mathbb{R}$ by
$\widetilde{f}(s)=\Lambda(f\cdot s)$ for each $f\in C(X)$. We
claim that $\widetilde{f}\in LUC(S)$. By Lemma~\ref{lemma}(ii),
the map $s\mapsto f\cdot s$ is norm continuous from $S$ to
$C(X)$. It is easy to verify that the continuity of $s\mapsto
\widetilde{f}\cdot s$ follows from the continuity of $s\mapsto
f\cdot s$.
 Define $J:C(X)\to \mathbb{R}$ by
$J(f)=I(\widetilde{f})$. Obviously $J$ is a left invariant
positive functional on $C(X)$. Therefore, $J(h)=0$ since $J$ is
invariant. But  $J(h)<0$ since $J$ is positive and $h<0$.

(ii)$\Rightarrow$(i):  It is easy to check that the set $M_U(S)$
of all means on $LUC(S)$  is a weak* compact subset of
$LUC(S)^*$.   Note that by Lemma~\ref{lemma}(i), the natural
action of $S$ from the left on $M_U(S)$ is continuous.  Let $\mu$
be a left $S$-invariant Radon probability measure  on $M_U(S)$.
Define $I_\mu:LUC(S)\to \mathbb{R}$ by $I_\mu(g)=\int \widehat{g}
d\mu$, where $\widehat{g}:M_U(S)\to \mathbb{R}$ is defined by
$\widehat{g}(J)=J(g)$. Clearly, $I_\mu$ is a left invariant mean
on $LUC(S)$.~\qed

\begin{rmk}
{\em If $S=G$ be locally compact group, then an invariant mean on
$LUC(G)$ extends to an invariant mean on the space $C_b(G)$ of
all bounded real-valued continuous functions on $G$ (cf.
\cite[Theorem~1.1.9, p.~21]{Runde}).}
\end{rmk}

A topological semigroup can be left, but not right, amenable
(e.g., consider the semigroup $S=\{a,b\}$ with the following
operation: $a\cdot a=b\cdot a=a,~a\cdot b=b\cdot b=b$). Of
course, if $S$ be a topological group, then $S$ is amenable if
and only if it is left (or right) amenable.  Basically it depends
on the fact that the operation $g\mapsto g^{-1}$ transposes the
order of products, and therefore interchanges left and right.
Also,  we will show that any abelian topological semigroup is
(both left and right) amenable (Corollary~\ref{abel-amen}).

Thanks to compactness of integral logic we have the following
 fact.

\begin{pro} \label{finite-amen}
Let $S$ be a topological semigroup with identity. Suppose that
there is a family $\{S_\alpha\}_{\alpha\in I}$ of subsemigroups of
 $S$ such that
\begin{itemize}

           \item [{\em (i)}] $\bigcup_{\alpha\in I}S_\alpha$ is dense in  $S$;
           \item [{\em (ii)}] $S_\alpha$ is an amenable subsemigroup with identity for all $\alpha\in I$;
           \item [{\em (iii)}]  For any $\alpha_1, \alpha_2\in I$, there exists $\alpha_3\in I$ such that $S_{\alpha_1}\bigcup S_{\alpha_2}\subseteq S_{\alpha_3}$.
\end{itemize}
  Then  $S$ is also amenable.
\end{pro}
{\bf Proof.} Let $S'=\bigcup_{\alpha\in I}S_\alpha$ and  $X$ be a
compact Hausdorff space and $\cdot$ a left   continuous action of
$S'$ on $X$. By assumptions, the theory $T_{S',X}$ of left
$S'$-invariant measures on $X$ is finitely satisfiable. By
Proposition \ref{Prop-1}, as $X$ and $\cdot$ are arbitrary, $S'$
is amenable. Assume that $I$ is an $S'$-invariant mean on
$LUC(S')$. Define $J:LUC(S)\to{\mathbb{R}}$ by
$J(f)=I(f\upharpoonright S')$ for each $f\in LUC(S)$. We can
easily check that $J$ is an left invariant mean on $LUC(S)$
because $S'$ is dense. Similarly, one can show that $S$ is right
amenable.~\qed

\begin{cor} \label{finite-amen2}
If every finitely generated subsemigroup (with identity) of a
topological semigroup $S$ is amenable, then $S$ is also amenable.
\end{cor}

Note that the converse may fail.  As an example let $S'$ be any
finitely generated non-amenable semigroup (e.g., the free group
on two generators), and let $S$ be a semigroup contains $S'$ and
one new element $s_0$ such that $s_0s=ss_0=s_0s_0=s_0$, for all
$s\in S'$. Then $S$ has an invariant mean $I(f)=f(s_0)$. The
subsemigroup $S'$ has not.

It is known that every locally compact group possesses a Haar
measure (cf.~\cite{Folland}), but not every locally compact group
is amenable. The free group on two generators, with the discrete
topology is a non-amenable locally compact group
(cf.~\cite[Example~449G, p.~399]{Fremlin4}). Of course, every
compact group is amenable. Indeed, assume that $G$ acts
continuously from the left on a compact Hausdorff space $X$. Fix
$x_0\in X$ and set $\phi(a) = a\cdot x_0$ for $a\in G$; then
$\phi$ is continuous. Let $\mu$ be the Haar probability measure
on $G$, and $\nu$ the Radon probability measure $\mu\phi^{-1}$ on
$X$. Clearly $\nu$ is $G$-invariant. As $X$ and $\cdot$ are
arbitrary, we have the following fact.

\begin{fct} \label{Haar}
Every compact group is amenable.
\end{fct}

A group $G$ is called {\em locally finite} if every finite subset
of $G$ generates a finite subgroup of $G$. An immediate
consequence of the above results is the following.

\begin{cor} \label{local-finit}
Let $G$ be a topological group such that the union of the finite
subsets of $G$ that generate a compact subgroup is dense.
  Then $G$ is amenable. In particular,
every locally finite topological group is amenable.
\end{cor}

\subsection{Commutativity}
The usual proof of the Bogolioubov-Krylov theorem uses the
Markov-Kakutani fixed point theorem. Now, we give a proof of this
theorem by using the compactness theorem and induction.

\begin{thm}[Bogolioubov-Krylov] \label{Krylov}
Assume that $S$ be an abelian semigroup which acts from the left
on a compact Hausdorff space $X$. Then
${\rm{Inv}}_X(S)\neq\emptyset$.
\end{thm}
{\bf Proof.} By Proposition \ref{Prop-1}, it suffices to consider
the case where $S$ is finite. We prove the theorem by induction
on the number of elements of $S$. Let $\cal D$ be a non-principal
ultrafilter on $\mathbb{N}$ and  $x_0$ any point of $X$. If
$S=\{s\}$, then define $\mu_1$ by $\int fd\mu_1=\lim_{n\to {\cal
D}} \frac{1}{n+1} \sum_{k=0}^n (f\cdot s^k)(x_0)$ for every $f\in
C(X)$. It is easy to check  that $\mu_1$ is
 invariant with respect to $s$. By induction hypothesis, there exists a
 measure $\nu$ on $X$ which is invariant with respect to
 $s_1,\ldots,s_{n-1}$.
 By the Riesz representation theorem, define the measure $\mu$ by $\int fd\mu=\lim_{n\to {\cal D}}
\frac{1}{n+1} \sum_{k=0}^n \int (f\cdot s_n^k)d\nu$ for every
$f\in C(X)$.  We can easily check that $\mu$ is invariant with
respect to $s_1,\ldots,s_n$. Indeed, it is easy to verify that
$\mu$ is $s_n$-invariant. Also, for each $i\leq n-1$, we have
\begin{align*}
\int (f\cdot s_i)d\mu &= \lim_{n\to {\cal D}} \frac{1}{n+1}
\sum_{k=0}^n \int (f\cdot s_i)\cdot
s_n^k d\nu\\
&= \lim_{n\to {\cal D}} \frac{1}{n+1} \sum_{k=0}^n \int (f\cdot
s_n^k)\cdot
s_i d\nu ~~~~~~~\text{commutativity}\\
&= \lim_{n\to {\cal D}} \frac{1}{n+1} \sum_{k=0}^n \int (f\cdot
s_n^k)d\nu ~~~~~~~~~~~\text{$\nu$ is $s_i$-invariant}\\&=\int
fd\mu.
\end{align*}
 Therefore, $\mu$ is the desired measure, so the theorem follows.~\qed

\bigskip
An immediate consequence of the Bogolioubov-Krylov theorem is the
following.

\begin{cor} \label{abel-amen}
 Any abelian topological  semigroup is amenable.
\end{cor}

Theorem \ref{Krylov} gives another proof of the existence of Haar
measure on abelian compact groups. By the same method  one can
also give a functional analytic proof of the existence of Haar
measures on abelian \textbf{locally} compact groups. We will
present a  proof of this theorem using the same method elsewhere.

\begin{cor}[Mazur-Orlicz]  Let ${\cal F}$ be a family of commuting mappings
of a set $X$ onto itself. Then there exists a mean on $B(X)$, the
space of all bounded real-valued functions on $X$,  which is $\cal
F$-invariant. In particular, every closed linear subspace $E$ of
$B(X)$ such that $f\circ h\in E$ whenever $f\in E$ and $h\in
{\cal F}$ has an $\cal F$-invariant mean.
\end{cor}
{\bf Proof.} Use Theorem \ref{Krylov}.~\qed

\subsection{Paradoxical decompositions}
The problematics of amenability has grown out of the famous
Banach-Tarski paradox (which essentially amounts to the
non-amenability of the free groups on two generators). We continue
this paper by looking at the connection between satisfiability
and paradoxical decompositions. Let $G$ be a discrete group
acting on a nonempty set $X$. Then $E\subseteq X$ is called {\em
$G$-paradoxical} if there are pairwise disjoint subsets
$A_1,\ldots,A_m,B_1,\ldots,B_n$ of $E$ along with $g_1,\ldots,
g_m,h_1,\ldots, h_n \in G$ such that $E =\bigcup_{i=1}^m g_i\cdot
A_i=\bigcup_{i=1}^n h_i\cdot B_i$. $X$ is said to be {\em
$G$-paradoxical} if it has a $G$-paradoxical subset. A group $G$
is called {\em paradoxical} if it is $G$-paradoxical. Clearly an
amenable group is non-paradoxical. A remarkable fact is that the
converse is also true, which follows from the following result of
Tarski.

\begin{thm} {\em  (\cite[p. 7]{Runde})} \label{T}
 Assume that $G$ and $X$ are as above. Then there exists a finitely additive, $G$-invariant
measure on $X$ defined for all subsets of $X$  if and only if $X$
is not $G$-paradoxical.
\end{thm}

A locally compact group $G$ admits a {\em Borel paradoxical
decomposition} if it has a paradoxical decomposition such that the
sets $A_1,\ldots,A_m,B_1,\ldots,B_n$ in the above definition are
Borel sets.  Paterson \cite{Paterson} proved that a locally
compact group $G$ is not amenable if and only if $G$ admits a
Borel paradoxical decomposition. The question of whether the
non-existence of such suitable paradoxical decompositions
characterizes the amenable, topological groups seems to be open
(cf. \cite{Wagon}).

 Now, we show that the amenability of a topological semigroup is expressible by a theory in integral logic.
 Note  that for a semigroup $S$ the dual  of the space $B(S)$  of all bounded real-valued functions on $S$  is the space of all signed charges
 on all subsets of $S$ (cf.~\cite[p. 496]{Ali}).
 Therefore, a mean $I$ on $B(S)$ is represented by a (positive)
charge $\nu_I$.  If $\nu_I$  is a charge which is not countably
additive, then $(S,\nu_I)$ is not a structure in integral logic.
Nevertheless, thanks to the representation theorem for
$M$-spaces, the amenability of a topological semigroup is
expressible. Indeed, consider a topological semigroup $S$  and
let $\sigma(S)$ ($=\sigma(LUC(S)))$ be the set of Riesz
homomorphisms $h:LUC(S)\to \mathbb{R}$ such that $h({\bf 1})=1$
(cf.~\cite[p.~222]{Fremlin3}). The set $\sigma(S)$  is sometimes
called the {\em spectrum} of $LUC(S)$. We will see that
$\sigma(S)$ is the space of complete types of a theory (see
Proposition~\ref{key} below). Note that, by~Proposition~353P(d) in
\cite[p.~243]{Fremlin3},
 $\sigma(S)$ is
the set $\texttt{M}_U(S)$ of all multiplicative means on $LUC(S)$.
 First, we remark that
$\sigma(S)$ is a weak* compact subset of $LUC(S)^*$ and $\|h\|=1$
for every $h\in\sigma(S)$, and hence by
 Lemma~\ref{lemma}(i), the natural action of $S$ on
$\sigma(S)$ is continuous.  The space $LUC(S)$ can be identified,
as normed Riesz space, with $C(\sigma(S))$, because $LUC(S)$ is
an $M$-space with standard order unit $\bf 1$ and $\sigma(S)$ is
a compact Hausdorff space (cf.~\cite[Corollary~354L]{Fremlin3}).
The identification is the map $f\mapsto \widehat{f}$ where
$\widehat{f}(h)=h(f)$ for $f\in LUC(S)$ and $h\in \sigma(S)$. By
the Riesz representation theorem, the identification of $LUC(S)$
with $C(\sigma(S))$ means that we have a one-to-one
correspondence $\mu\leftrightarrow I_\mu$ between Radon
probability measures $\mu$ on $\sigma(S)$ and positive linear
functionals $I_\mu$ on $LUC(S)$ such that $I_\mu({\bf 1}) = 1$,
given by the formula $I_\mu(f) = \int \widehat{f} d\mu$ for $f\in
LUC(S)$. Now
\begin{align*}
\mbox{$I_\mu$ is invariant} &\Leftrightarrow I_\mu(f\cdot s)=I_\mu(f)~\mbox{ for every $f\in LUC(S)$ and $s\in S$}\\
&\Leftrightarrow  \int \widehat{f\cdot s}d\mu=\int \widehat{f} d\mu~\mbox{ for every $f\in LUC(S)$ and $s\in S$} \\
&\Leftrightarrow \int (\widehat{f}\cdot s)d\mu=\int \widehat{f}
d\mu~\mbox{ for every $f\in LUC(S)$ and $s\in S$}
\\&\Leftrightarrow \mbox{$\mu$ is invariant}.
\end{align*}
So there is a one-to-one correspondence between Radon probability
left $S$-invariant measures on $\sigma(S)$ and left $S$-invariant
means on $LUC(S)$. Let  $T_S=T_{S,\sigma(S)}$  be the theory of
left $S$-invariant measures on $\sigma(S)$. Summarizing, we have
the following.

\begin{pro} \label{amen-paradox-satis}
Assume that $S$ and $T_S$ are as above. Then the following are
equivalent:
\begin{itemize}
           \item [{\em (i)}]  $S$ is  amenable.
           \item [{\em (ii)}] $T_S$ is  satisfiable.
\end{itemize}
         If   $S$ is a locally compact group, then {\em (i)}
and {\em (ii)} are equivalent to
\begin{itemize}
           \item [{\em (iii)}]  $S$ is not Borel paradoxical.
\end{itemize}
\end{pro}

In fact  we can say more:  if $S$ and $T_S$ are as above, then the
cardinal of the set of all left $S$-invariant means on $LUC(S)$
is equal to the number of models of $T_S$ up to almost
isomorphism. Indeed, if $\mu\neq\nu$ are (left) $S$-invariant
measures on $\sigma(S)$ then $(\sigma(S),{\cal B}, \mu)$ and
$(\sigma(S),{\cal B}, \nu)$ with the natural interpretation of
relation and constant symbols are different models of $T_S$.
Conversely, assume that ${\M}=(M,{\cal B},\mu_{\M})$ is a model
of $T_S$.  By Proposition \ref{Pushing-Down}, the substructure
${\M}'=(\sigma(S),{\cal B}_{\sigma(S)},
\mu_{\M}\upharpoonright{\sigma(S)})$ is also a model of $T_S$ and
the inclusion map $\sigma(S)\to M$ covers a full measure subset
of $M$. Therefore, ${\M}'\simeq_a{\M}$. Clearly, the unique
extension of $\mu_{\M}\upharpoonright \sigma(S)$ to a Radon
measure on $\sigma(S)$ is left $S$-invariant.  To summarize:

\begin{pro} \label{cardinal}
Assume that $S$ and $T_S$ are as above. Then there
 is a bijection between the set of all models of $T_S$ and the set of all left $S$-invariant means on $LUC(S)$.
\end{pro}

\subsection{Extreme amenability} \label{5}
In this subsection we present some other results  for extremely
amenable  topological semigroups.  Most of the proofs are
straightforward and we omit some unnecessary details.
 First, we  characterize extremely  amenable topological
semigroups in terms of  multiplicative invariant measures
(Fact~\ref{2}). Finally, we prove that the extreme amenability of
a topological semigroup is expressible by a theory in integral
logic (Proposition \ref{extr-amen}).

A Radon probability measure $\mu$ on a compact Hausdorff space $X$
is {\em multiplicative} if $\int fd\mu\times \int gd\mu=\int
(f\times g)d\mu$ (the pointwise product) for all $f,g\in C(X)$.

Let $S$ be a topological semigroup which acts on a compact
hausdorff space $X$ from the left. Let $\texttt{T}_{S,X}$ be the
theory of left $S$-invariant measures on $X$ with the additional
axiom schema
\begin{itemize}
\item [(9)] $ \int R_{f\times g}(x)dx=\int R_{f}(x)dx\times\int R_{g}(x)dx$ \ \ \ \ \ \   ~ for each $R_f,R_g,R_{f\times g}\in L_{X}$,

 where  $(f\times g)(x)=f(x)\times g(x)$.
\end{itemize}
 \vspace{.3 cm}
 Note that (9)  says that the measure is multiplicative.
 $\texttt{T}_{S,X}$ is called the theory of {\em
multiplicative left $S$-invariant measures on $X$}.

Let ${\rm{MInv}}_X(S)$ be the set of all multiplicative, Radon
probability measures on X which are left $S$-invariant. A
consequence of the compactness theorem is the following.

\begin{pro} \label{Prop-1-ext}
Assume that $S, X$ and $\texttt{T}_{S,X}$ are as above. Then the
following are equivalent:
\begin{itemize}
           \item [{\em (i)}] ${\rm{MInv}}_X(S)\neq\emptyset$.
           \item [{\em(ii)}] $\texttt{T}_{S,X}$ is  satisfiable.
\end{itemize}
\end{pro}

Let $S$ be a topological semigroup. A mean $I$ on $LUC(S)$ is {\em
multiplicative} if $I(f)\times I(g) = I(f\times g)$ (the
pointwise product) for all $f,g\in LUC(S)$. We remark  that
$LUC(S)$ is a closed and invariant subalgebra of $C_b(S)$
(cf.~\cite[Lemmas~1.1 and~1.2]{Namioka}).

\begin{dfn} \label{ext-semi-amen}
 {\em    A topological semigroup $S$ is said to be {\em
extremely left} {\em (right)} {\em amenable} if $LUC(S)$
($RUC(S)$) admits a multiplicative left (right) invariant mean. A
topological semigroup $S$ is called {\em extremely amenable} if
it is both left and right amenable.}
\end{dfn}

 \begin{rmk}
 {\em  A topological semigroup  $S$ has the {\em left (right) fixed
point on compacta property} if every continuous action of $S$ on
a compact Hausdorff space by the left (right) side has a fixed
point. In \cite{Mitch2}, Mitchell showed  that a topological
semigroup $S$ has a multiplicative left invariant mean on
$LUC(S)$ iff $S$ has the left fixed point on compacta property.
Also, he asked the question: Is there a non trivial extremely
amenable group at all? Historically the first example of extremely
amenable groups was found in \cite{Herer}. Many further examples
of extremely amenable groups may be found in
\cite{Pestov1,Pestov2,Fremlin4}.}
\end{rmk}

The following fact presents a proof of Mitchell's theorem
\cite[Theorem~1]{Mitch2} and it also characterizes extremely
amenable topological semigroups in terms of multiplicative
invariant measures.

\begin{fct} \label{2}
Let $S$ be a topological semigroup with identity. Then the
following are equivalent:
\begin{itemize}
             \item [{\em (i)}]  $S$  is extremely left amenable.
             \item [{\em (ii)}] $S$ has the left  fixed point on compacta property.
             \item [{\em (iii)}] Whenever  $X$ is a
non-empty compact Hausdorff space and $\cdot$ is a  continuous
action of $S$ on $X$  by the left   side, then
${\rm{MInv}}_X(S)\neq\emptyset$.
\end{itemize}
\end{fct}
{\bf Proof.}  (i)$\Leftrightarrow$(iii): The set
$\texttt{M}_U(S)$ ($=\sigma(S)$) of all multiplicative means on
$LUC(S)$ is a weak* compact subset of $LUC(S)^*$. By
Lemma~\ref{lemma}(i), the natural action of $S$ on
$\texttt{M}_U(S)$ (by the left side) is continuous. Also, it is
easy to verify that ${\rm{MInv}}_X(S)\neq\emptyset$ iff for every
elements $s_1,\ldots,s_n$ of $S$ and elements
$f_1,\ldots,f_n,g_1,\ldots,g_n$ of $C(X)$ we have $\| {\bf
1}-\sum_{i=1}^n g_i\times(f_i\cdot s_i-f_i)\|\geq 1$. (Compare
Fact~\ref{Compact-Haar-3}.)  Now, the proof is a simple adaptation
of the proof of Fact~\ref{amenable-equiv}.

(ii)$\Rightarrow$(iii): Assume that $x_0\in X$ is a fixed point,
i.e., $s\cdot x_0=x_0$ for every $s\in S$. Define the measure
$\mu$ by $\int fd\mu=f(x_0)$ for every $f\in C(X)$. Clearly,
$\mu$ is a multiplicative Radon left $S$-invariant measure on $X$.

(iii)$\Rightarrow$(ii): Assume that $X$ is a non-empty compact
Hausdorff space and $\cdot$ is a  continuous action of $S$ on
$X$  by the left side. Let $\mu$ be a multiplicative left
$S$-invariant Radon probability measure on $X$. Then the linear
functional $I$ defined by $I(f)=\int fd\mu$ is multiplicative and
invariant. Therefore, by Lemma~25 in \cite[p.~278]{Dunford},
there is a point $x_0$ in  $X$ such that $I(f)=f(x_0)$ for every
$f\in C(X)$. Since $C(X)$ separates points and $I$ is invariant,
$x_0$ is the desired fixed point.~\qed

\bigskip
Using the compactness theorem of integral logic, one can prove the
following fact.

\begin{pro} \label{ext-finite-amen}
 If $S$ is a
topological semigroup with a dense subset $\bigcup_{\alpha\in
I}S_\alpha$ where $S_\alpha$ are extremely amenable
  semigroups and for any $\alpha_1,\alpha_2\in I$, $S_{\alpha_1}\bigcup S_{\alpha_2}\subseteq S_{\alpha_3}$ for some $\alpha_3\in I$ then
 $S$ is extremely amenable.
\end{pro}

At the end of this section we show that the extreme amenability of
a topological semigroup is expressible by a theory in integral
logic. Let $S$ be a topological  semigroup and
$\texttt{T}_{S}=\texttt{T}_{S,\sigma(S)}$ be the theory of
multiplicative left $S$-invariant measures on $\sigma(S)$. In
fact, we show that the cardinal of ${\rm{MInv}}_X(S)$ is  equal
to the number of models of $\texttt{T}_{S}$. By
Propositions~\ref{amen-paradox-satis} and \ref{cardinal}, it
suffices to show that  there is a one-to-one correspondence
  between  multiplicative Radon probability measures on
$\sigma(S)$ and multiplicative means on $LUC(S)$. Note that the
identification of $LUC(S)$ and $C(\sigma(S))$ is algebraic, i.e.,
$\widehat{f\times g}=\widehat{f}\times \widehat{g}$ for every
$f,g\in LUC(S)$ (cf.~\cite[Pro~353P(d), p.~243]{Fremlin3}).
 Now
\begin{align*}
\mbox{$I_\mu$ is multiplicative} &\Leftrightarrow I_\mu(f\times g)=I_\mu(f)\times I_\mu(g)~\mbox{ for every $f,g\in LUC(S)$}\\
&\Leftrightarrow  \int \widehat{f\times
g}d\mu=\int \widehat{f} d\mu\times\int \widehat{g} d\mu~\mbox{ for every $f,g\in LUC(S)$} \\
&\Leftrightarrow \int (\widehat{f}\times \widehat{g})d\mu=\int
\widehat{f} d\mu\times\int \widehat{g} d\mu~\mbox{ for every
$f,g\in LUC(S)$}
\\&\Leftrightarrow \mbox{$\mu$ is multiplicative}.
\end{align*}
 To summarize:

\begin{pro} \label{extr-amen}
Assume that $S$ and $\texttt{T}_{S}$ are as above. Then there is
a bijection between the set of all models of $\texttt{T}_{S}$ and
the set of all multiplicative left $S$-invariant means on
$LUC(S)$. In particular,  $S$ is extremely left amenable iff
$\texttt{T}_{S}$ is  satisfiable.
\end{pro}

\section{Types and stability} \label{6}
In classical model theory, a complete type determines a finitely
additive $0$-$1$ valued  measure on the formulas. Actually, one
can say more, i.e., a complete type is a $0$-$1$ valued  Riesz
homomorphism on the formulas. Indeed, let $L$ be a first order
language,  $\M$ an $L$-structure, $a$ an element of $M$, and
$\text{tp}^{\M}(a)$ be the complete type of $a$ in $\M$. For each
$L$-formula $\phi(x)$, define $f_\phi:M\to\{0,1\}$ by
$f_\phi(b)=1$ if ${\M}\vDash \phi(b)$, and $f_\phi(b)=0$
otherwise. Let $V=\{f_\phi: \phi\in L\}$. One can easily check
that $V$ is an (Archimedean) Riesz space (see
Definitions~\ref{Riesz space} and \ref{Riesz homomorphism} below).
For this we define $f_\phi+f_\psi:=f_{\phi\vee\psi}$,
$-f_\phi:=f_{\neg\phi}$, and for each $r\in \mathbb{R}$, $r\cdot
f_{\phi}:=f_\phi$ if $r>0$, $r\cdot f_{\phi}:=f_{\neg\phi}$ if
$r<0$, and  $r\cdot f_{\phi}:=\bf{0}$ if $r=0$. Also, $f_\phi\leq
f_\psi$ if $f_\phi(b)\leq f_\psi(b)$ for each $b\in M$. Clearly,
$V$ with this structure is a Riesz space, i.e.,  it is a partially
ordered linear space which is a lattice.
 Now, for an $a\in M$, define the Riesz homomorphism
$I_a:V\to\{0,1\}$ by $I_a(f_\phi)=1$ if $f_\phi(a)=1$, and
$I_a(f_\phi)=0$ otherwise, i.e., $I_a(f_\phi)=1$ iff
$\phi\in\text{tp}^{\M}(a)$. In other words, $I_a$ can be
interpreted as playing the role of $\text{tp}^{\M}(a)$.

More generally, we consider real valued Riesz homomorphisms.
Indeed, consider an arbitrary partially ordered set ${\cal
L}=\{f_\phi:M\to {\mathbb R}: \phi\in L\}$ such that
\begin{align*}
\forall  b\in M~ \colon ~f_\phi(b)\leq f_\psi(b) ~ &\Longleftrightarrow ~  \vDash\phi(b)\to\psi(b)\\
  ~f_\phi(b)< f_\psi(b) ~ &\Longleftrightarrow ~
\vDash\neg\phi(b)\wedge\psi(b).
\end{align*}
Let $V$ be the linear space generated by $\cal L$. Again, $V$ is
an Archimedean Riesz space. Define the Riesz homomorphism
$I_a:V\to\mathbb R $ by $I_a(f)=f(a)$. It is easy to verify that
$\phi\in\text{tp}^{\M}(a)$ iff $I_a(f_{\phi\vee\neg\phi})\leq
I_a(f_\phi)$. Therefore  it is natural to conjecture that real
valued Riesz homomorphisms on measurable functions should play
the role of complete types in the framework of integral logic. Our
next goal is to convince the reader that this is indeed the case.

\subsection{Types}   Let us now return to integral logic.
Suppose that $L$ is an arbitrary language, maybe with $n$-ary
relation symbols and $n$-ary function symbols. Let $\M$ be a {\em
graded} $L$-structure as discussed in \cite{BP}, $A\subseteq M$
and $T_A=Th({\M}, a)_{a\in A}$. Let $p(x)$ be a set of
$L(A)$-statements in free variable $x$. We shall say that $p(x)$
is a {\em type  over} $A$ if $p(x)\cup T_A$ is satisfiable. A
{\em complete type over} $A$ is a maximal type over $A$.  We let
$S^{\M}(A)$ be the set of all complete types over $A$. The {\em
type of $a$ in $M$ over $A$}, denoted by $\text{tp}^{\M}(a/A)$,
is the set of all $L(A)$-statements satisfied in $\M$ by $a$. For
$\phi(x)$ an $L(A)$-formula, we let $$[\phi>0]=\{p\in S^{\M}(A):
\mbox{ for some $\epsilon>0$ the statement $(\phi\geq\epsilon)$ is
in $p$ }\}.$$
 The {\em logic topology} (or the
{\em Stone topology}) on $S^{\M}(A)$ is the topology generated by
taking the sets $[\phi>0]$ as basic open sets.
 We will give a characterization of the complete types. First, we need some
notions from functional analysis.

\begin{dfn}[Riesz space] \label{Riesz space}
{\em A  {\em Riesz space} or {\em vector lattice} is a partially
ordered linear space which is a lattice. A Riesz space
$\mathcal{L}$ is called {\em Archimedean} if
$\inf_{\delta>0}\delta f=\textbf{0}$ for each $f\geq \textbf{0}$
in $\mathcal{L}$. An element $\textbf{1}\geq\textbf{0}$ of
$\mathcal{L}$ is an {\em order unit}  in $\mathcal{L}$ if for
every $f\in\mathcal{L}$ there is an $n\in \mathbb{N}$ such that
$|f|\leq n \textbf{1}$. }
\end{dfn}

The following notion will play a fundamental role in what follows.

\begin{dfn}[Riesz homomorphism] \label{Riesz homomorphism}
{\em Let $\mathcal{L}, \mathcal{L}'$ be partially ordered linear spaces.
 A {\em Riesz homomorphism} from $\mathcal{L}$ to $\mathcal{L}'$ is
 a linear operator $T : \mathcal{L}\to \mathcal{L}'$ such that whenever $A
\subset \mathcal{L}$ is a finite non-empty set and $\inf A
=\textbf{0}$ in $\mathcal{L}$, then $\inf T[A] = \textbf{0}$ in
$\mathcal{L}'$. }
\end{dfn}

Any Riesz homomorphism is a {\em positive} linear operator, i.e.
$T(f)\geq \textbf{0}$ for all $f\geq \textbf{0}$ (see
\cite{Fremlin3}, 351H(b)).

\begin{fct}[\cite{Fremlin3}, 354K] \label{embed to C(X)}
Let $\mathcal{L}$ be an Archimedean Riesz space with order unit $\textbf{1}$.
 Then it can be embedded as an order-dense and norm-dense Riesz subspace of $C(X)$, where $X$ is
a compact Hausdorff space, in such a way that $\textbf{1}$
corresponds to $\chi_X$; moreover, this embedding is essentially
unique.
\end{fct}

The compact space $X$ in Fact~\ref{embed to C(X)} is the set of
Riesz homomorphisms $I$ from $\mathcal{L}$ to $\mathbb{R}$ such
that $I(\textbf{1})=1$, and the embedding is the map $T :{\cal
L}\to \mathbb{R}^X$ defined by setting $(Tf)(I) = I(f)$ for every
$I\in X$, $f\in\mathcal{L}$ (see the proof of Theorem~353M in
\cite{Fremlin3}).

\medskip
Let $\M$ be an $L$-structure and $A$ a subset of $M$. We define
  ${\cal L}_A$ to be the family of all measurable functions $\phi^{\M}$
 where $\phi$ is an $L(A)$-formula with a free variable
$x$ (see  the paragraph after Definition~\ref{structure}). Then
${\cal L}_A$ has a natural Riesz space structure given by
$(\phi^{\M}+\psi^{\M})(a)=\phi^{\M}(a)+\psi^{\M}(a)$,
$(r\phi^{\M})(a)=r\phi^{\M}(a)$ for all $a\in M$, and
$\phi^{\M}\geq\psi^{\M}$ iff $\phi^{\M}(a)\geq\psi^{\M}(a)$  for
all $a\in M$. Also, $|\phi^{\M}|(a)=|\phi^{\M}(a)|$ for all $a\in
M$,  $\min(\phi^{\M},\psi^{\M})$, $\max(\phi^{\M},\psi^{\M})$ are
in ${\cal L}_A$, and $\|\phi^{\M}\|=\sup_{a\in M}|\phi^{\M}(a)|$.
Clearly, ${\cal L}_A$ is Archimedean. The constant function
$\textbf{1}$ is an order unit and the uniform norm is its
order-unit norm (see \cite[354G(a)]{Fremlin3}).

Let $\sigma_A({\M})$ be the set of Riesz homomorphisms $I: {\cal
L}_A\to \mathbb{R}$ such that $I(\textbf{1}) = 1$. This set
 is called the {\em spectrum} of $T_A$. Since
${\cal L}_A$ is a normed linear space (with the uniform norm),
the unit ball $B^*=\{I\in{\cal L}_A^*:\|I\|\leq 1\}$ in ${\cal
L}_A^*$  is compact in the weak* topology by Alaoglu's Theorem.
Also, we know that $\sigma_A({\M})$ is the set of positive {\em
extreme points} of the unit ball  $B^*$, i.e.
$\sigma_A({\M})=\{I\in B^*:\|I\|=1 \mbox{ and  $I$ is
positive}\}$ (see \cite{Fremlin3}, 354Y(j)). Since
$\sigma_A({\M})\subseteq B^*$ is weak* closed, so it is weak*
compact. (We remark that the weak* topology on $\sigma_A({\M})$ is
simply the topology of pointwise convergence: $I_\alpha\to I$ in
the weak* topology iff $I_\alpha(\phi^{\M})\to I(\phi^{\M})$ for
all $\phi^{\M}\in {\cal L}_A$; see \cite{Folland}, page 169, for
details.)



\medskip
 The next propositions show that a complete
type can be coded by a Riesz homomorphism and give a
characterization of complete types. The key idea behind these
propositions is a construction which allows us to consider $\M$ as
an elementary submodel of the type space $S^{\M}(M)$ with the
appropriate structure.

\begin{dfn}[$\sigma_M({\M})$ as an elementary extension]  \label{construction} {\em Assume that $\M$ is an $L$-structure and $\mu$ is the
measure on $M$. By Fact~\ref{embed to C(X)}, the space ${\cal
L}_M$ can be embedded as an order-dense and norm-dense Riesz
subspace of $C(\sigma_M(\M))$.  The embedding is the map $T :{\cal
L}_M\to \mathbb{R}^{\sigma_M({\M})}$ defined by setting
$(T\phi^{\M})(I) = I(\phi^{\M})$ for every $I\in
\sigma_M({\M})$,  $\phi^{\M}\in\mathcal{L}_M$. We define the
elementary extension ${\N}=(\sigma_M({\M}),\nu,
T\phi^{\M})_{\phi^{\M}\in\mathcal{L}_M}$ of $\M$ with the natural
interpretations of symbols and measure as follows:

First, we can easily see that $M\subseteq \sigma_M({\M})$.
Indeed, for each $a\in M$, define
$I_a:\mathcal{L}_M\to\mathbb{R}$ by $I_a(\phi^{\M})=\phi^{\M}(a)$
for $\phi^{\M}\in\mathcal{L}_M$. Now, one can assume that  the
language has a $2$-ary relation symbol $\textbf{e}$  with the
interpretation $\textbf{e}(a,b)=1$ if $a=b$, and
$\textbf{e}(a,b)=0$ otherwise (cf. \cite[p. 469]{BP}). Therefore,
$I_a\neq I_b$ if $a\neq b \in M$. More generally, if ${\cal L}_M$
separates $M$, i.e. for each $a\neq b\in M$ there is
$\phi^{\M}\in{\cal L}_M$ such that
$\phi^{\M}(a)\neq\phi^{\M}(b)$, then $I_a\neq I_b$. To summarize,
the map $M\hookrightarrow \sigma_M({\M})$ defined by $a\mapsto
I_a$ is injective, and so we can assume that $a=I_a$ and
$M\subseteq\sigma_M({\M})$.

Second, define $\nu\{T\phi^{\M}>0\}:=\mu\{\phi^{\M}>0\}$ for all
$\phi^{\M}\in\mathcal{L}_M$. (Recall that $\{T\phi^{\M}>0\}$ is
the set $\{I\in \sigma_M({\M}): T\phi^{\M}(I)>0\}$.) Then $\nu$ is
a premeasure on the algebra ${\cal A}=\{\{T\phi^{\M}>0\}:
\phi^{\M}\in\mathcal{L}_M\}$. By Carath\'{e}odory's theorem,
$\nu$ has a unique extension to a measure on the $\sigma$-algebra
generated by $\cal A$, still denoted by $\nu$. Also, we can
assume that $M$ is $\nu$-measurable and $\nu(N\setminus M)=0$,
i.e. $M$ has full-measure.

Third, for each formula $\phi(x,y_1,\ldots,y_n)$ and elements
$a_1,\ldots,a_n$ of $M$,
 define $\phi^{\N}(x,I_{a_1},\ldots,I_{a_n}):\sigma_M({\M})\to\mathbb{R}$ by
  $\phi^{\N}(x,I_{a_1},\ldots,I_{a_n})=T\phi^{\M}(x,a_1,\ldots,a_n)$.
  Then for each $b\in M$ we have
\begin{align*}
   \phi^{\N}(I_b,I_{a_1},\ldots,I_{a_n}) &
   =T\phi^{\M}(x,a_1,\ldots,a_n)(I_b)
   \\ & =I_b(\phi^{\M}(x,a_1,\ldots,a_n)) \\ &
       =\phi^{\M}(b,a_1,\ldots,a_n).
\end{align*}
Also, for a formula $\phi(x_1,x_2)$, define
$\phi^{\N}(x_1,x_2):(\sigma_M({\M}))^2\to\mathbb{R}$ by
$\phi^{\N}(I_a,I)=T\phi^{\M}(a,y)(I)$ and
$\phi^{\N}(I,I_b)=T\phi^{\M}(x,b)(I)$, where $a,b\in M$ and $I\in
N$, and $\phi^{\N}(I,I')=0$ if $I,I'\in N\setminus M$. Similarly,
we can define $\phi^{\N}(x_1,\ldots,x_n)$. For a $2$-ary function
symbol $f$, define $f^{\N}(I_a,I_b):=f^{\M}(a,b)$ for all $a,b\in
M$, and  for some $I''\in N\setminus M$, $f^{\N}(I,I'):=I''$ if
at least one of $I,I'$ belongs to $N\setminus M$.
 Similarly, we can define $f^{\N}(x_1,\ldots,x_n)$. Also, we can assume that the $n$-ary
relations and functions on $N$ are $\nu_n$-measurable. In fact,
our definitions are not important on the set $N^n\setminus M^n$,
because $\nu_n(N^n\setminus M^n)=0$ and we can take an
appropriate $\sigma$-algebra on $N^n$. }
\end{dfn}

\begin{pro}
Assume that $\M$, $\N$ are as above. Then ${\M}\preceq{\N}$.
\end{pro}
{\bf Proof.} Since $M\subseteq \sigma_M({\M})$ and
$\phi^{\N}(\b)=\phi^{\M}(\b)$ for all $\b\in M$ and formula
$\phi(\x)$, so $\M$ is a substructure of $\N$. Now by the
Tarski-Vaught test (Proposition~\ref{Tarski-Vaught} above), $\N$
is an elementary extension of $\M$. Indeed, we note that
$\nu\{\phi^{\N}>0\}=\mu\{\phi^{\M}>0\}$ for all
$\phi^{\M}\in\mathcal{L}_M$. (See also \cite{BP},
Proposition~5.10.)~\qed

\bigskip
Also, we will see that $\N$ realizes every type in
$S^{\M}(M)$; in fact $S^{\M}(M)=\sigma_M({\M})$.

\begin{pro} \label{key}
Assume that $\M$ is an $L$-structure and $A\subseteq M$.
\begin{itemize}
             \item [{\em (i)}]  There is a bijection from $S^{\M}(M)$ onto $\sigma_M({\M})$.
             \item [{\em (ii)}] $q\in S^{\M}(A)$ if and only if there is
             an elementary extension $\N$ of $\M$ and
              $x_0\in N$ such that $q=\text{tp}^{\N}(x_0/A)$.
\end{itemize}
\end{pro}
{\bf Proof.}  (i):  Assume that $p(x)$ is a complete type over
$\M$. Define $I_p: {\cal L}_M\to \mathbb{R}$ by $I_p(\phi^{\M})=r$
if the statement $\phi(x)=r$ is in $p(x)$. Clearly, $I_p$ is a
Riesz homomorphism on ${\cal L}_M$ and $I_p(\textbf{1})=1$. The
map $p\mapsto I_p$ is injective, and we may reasonably assume
that $p=I_p\in \sigma_M({\M})$. In particular, for any $a\in M$,
$\text{tp}^{\M}(a/M)=\{\phi(x)=\phi^{\M}(a): \phi\in {\cal L}_M\}$
and $I_{\text{tp}^{\M}(a/M)}(\phi^{\M})=\phi^{\M}(a)$. (Before we
showed that the map $M\hookrightarrow \sigma_M({\M})$ defined by
$a\mapsto I_{\text{tp}^{\M}(a/M)}$ is  injective.)

Now, we show that the map $p\mapsto I_p$ is surjective. Assume
that $I\in\sigma_M(\M)$. Let ${\N}=(\sigma_M({\M}),\nu,
T\phi^{\M})_{\phi^{\M}\in\mathcal{L}_M}$ be the elementary
extension of $\M$ constructed in Definition~\ref{construction}
and $p=\text{tp}^{\N}(I/M)$. Then, it is easy to check that
$I_p=I$. (Indeed, recall that
$\phi^{\N}(I)=T\phi^{\M}(I)=I(\phi^{\M})$ for all
$\phi^{\M}\in{\cal L}_M$.) Therefore, the map $p\mapsto I_p$ is
also surjective.

(ii): Let $q\in S^{\M}(A)$ and $\N$ be the elementary extension of
$\M$ constructed in Definition~\ref{construction}.  Assume that
$p\in S^{\N}(M)=S^{\M}(M)$ is an extension of $q$. Then there is a
point $x_0\in N$ such that $p=\text{tp}^{\N}(x_0/M)$ (see (i)
above). Clearly, $q=\text{tp}^{\N}(x_0/A)$.~\qed

\bigskip
Recall that $\sigma_M({\M})$   is weak* compact. Since
$S^{\M}(M)=\sigma_M({\M})$, we can also equip $S^{\M}(M)$ with the
weak* topology. It is easy to check that the weak* topology and
the logic topology on $S^{\M}(M)$ are the same. Indeed, for each
$\phi^{\M}\in {\cal L}_M$, define
$\phi:S^{\M}(M)\to[-\flat_\phi,\flat_\phi]$ by $p\mapsto
I_p(\phi^{\M})$. Then obviously the logic topology on $S^{\M}(M)$
is the weakest topology in which all the functions $p\mapsto
\phi(p)$ are continuous. Therefore, for arbitrary $p_\alpha,p\in
S^{\M}(M)$
\begin{align*}
    I_{p_\alpha}\to I_p \mbox{ in the weak* topology } \ \  &
   \Leftrightarrow \ \  I_{p_\alpha}(\phi^{\M})\to I_p(\phi^{\M})
   \mbox{ for all } \phi^{\M}\in{\cal L}_M
   \\ & \Leftrightarrow  \ \   \phi(p_{\alpha})\to \phi(p) \mbox{ for all } \phi^{\M}\in{\cal
   L}_M  \\ &
    \Leftrightarrow  \ \  \phi \mbox{ is continuous, for all } \phi^{\M}\in{\cal L}_M
       \\ & \Leftrightarrow \ \  p_\alpha\to p  \mbox{ in the logic topology. }
\end{align*}

\begin{rmk} \label{topomeasure}
{\em  By Proposition~\ref{key}, the elementary extension
${\N}=(\sigma_M({\M}),\nu,\phi^{\N})$, as constructed in
Definition~\ref{construction}, realizes every type over $M$. Also,
it is easy to verify that $M$ is a dense subset of
$N=\sigma_M({\M})$. Indeed, if $M$ is not dense in $N$, there is
a non-zero $h\in C(N)$ such that $h(I_a) = 0$ for every $a\in M$;
but as the uniform completion $\overline{{\cal L}}_M$ of ${\cal
L}_M$ is identified with $C(N)$ (because ${\cal L}_M$ is dense in
$C(N)$), there is an $f \in \overline{{\cal L}}_M$ such that
$I(f) = h(I)$ for every $I\in N$. Assume that $f_n\to f$
uniformly, where $f_n\in {\cal L}_M$. Therefore, there are
$h_n\in C(N)$ such that $I(f_n) = h_n(I)$ for every $I\in N$.
Clearly, $h_n\to h$ uniformly. In this case, $f$ cannot be the
zero function, but $f(a)=\lim_n f_n(a)=\lim_n I_a(f_n)=\lim_n
h_n(I_a)=h(I_a)=0$ for every $a\in M$. Thus the image of $M$ is
dense, as claimed.

 On the other hand, since $\phi^{\N}$'s are continuous,
  the natural  measure $\nu$ on $\N$ is Baire  and  it has a
unique extension to a Radon measure, which again we denote this
measure by $\mu$. From now on we assume that ${\N}=(S({\M}),\mu)$
with the appropriate structure, where $\mu$ is this Radon
measure.}
\end{rmk}

\begin{cor}
Let $G$ be an amenable topological group and $T_G$ the theory of
left $G$-invariant  measures on $\sigma(G)$. Then $G$ is
extremely amenable iff there is a complete type $p\in
S(\sigma(G))$ such that $g\cdot p=p$ for each $g\in G$.
\end{cor}

\subsection{Definable relations}
\begin{dfn} \label{Definable relations}{\em A
relation $\xi:M\rightarrow[-\flat,\flat]$ is {\em
$\emptyset$-definable} if there is a sequence $\phi_k(\x)$ of
formulas such that $\flat_{\phi_k}\leq\flat$ and
$\phi_k\rightarrow \xi$ pointwise. A subset is definable if its
characteristic function is definable.}
\end{dfn}

This may be defined on the basis of other notions of convergence
such as almost uniform convergence, convergence in measure,
convergence in the mean etc. However, the corresponding
definitions are equivalent. For example if $\phi_k$ converges in
measure to $\xi$, then it has a subsequence which converges to
$f$ almost everywhere. So, if $R$ is definable using the first
notion of convergence, it is also definable using the second one.
In particular, since the measure is finite and $|\phi_k|\leq
\flat$, $\phi_k\rightarrow \xi$ in measure iff $\phi_k\rightarrow
\xi$ in mean iff $\phi_k\rightarrow \xi$ pointwise (see
\cite{Folland}).  On the other hand, if $\M\preceq\N$ and $\xi$ is
definable in $\M$, then there is a corresponding definable
relation $\xi'$ in $\N$ and it is not hard to see that
${\M}\preceq_{a}\N$. The set of definable relations is a Banach
algebra with the norm defined by $\|\phi\|=\sup_{x}|\phi(x)|$ and
this algebra depends only on $T$. It can be described as the
completion of the algebra of formulas with the  uniform norm. We
denote this completion by $L(T)$.  A relation is $\M$-definable
if it is definable in $Th({\M},a)_{a\in M}$. So, $L(\M)$ is
defined in the natural way.

\subsection{Local stability}
 Here and in the next section we give two different notions of ``stability" of a formula inside a
model, a measure theoretic notion and a model theoretic notion.
In fact, the measure theoretic notion (Definition~\ref{stab1}) is
a suitable form of the dependence property in classical model
theory.

Let $\M$ be a structure and $\phi(x,y)$ a formula. Assume that
${\N}\succeq{\M}$ and $a\in N$. Let $p =\text{tp}_\phi^{\M}(a/M)$
be the {\em complete $\phi$-type of $a$ over $M$}, i.e., a
function which associates to each instance $\phi(x,b)$, $b\in M$,
the value $\phi(a,b)$, which will then be denoted by $\phi(p,b)$.
Note that the complete $\phi$-type $p$ uniquely determines a Riesz
homomorphism $I_{p}:{\cal L}_{\phi}\to \mathbb{R}$ where ${\cal
L}_\phi$ is the Riesz space generated by $\{\phi(x,b):b\in M\}$,
and $I_p(\phi(x,b))=\phi(p,b)$ for each $b\in M$.
 We equip
$S_\phi(M)$ with the weakest topology in which all functions
$p\mapsto \phi(p,b),b\in M$ are continuous. Equivalently, if
$\sigma_\phi(M)$ be the spectrum of $T_\phi=\{\phi\geq r:\phi\geq
r \mbox{ is in }T({\M},a)_{a\in M}\}$ (i.e., the set of Riesz
homomorphisms $I : {\cal L}_\phi \to \mathbb{R}$ such that
$I({\bf 1}) = 1$), then $S_\phi(M)=\sigma_\phi(M)$ is equipped
with the topology induced by the weak* topology on ${\cal
L}^*_\phi$.
 Clearly, $S_\phi(M)$ is
compact Hausdorff. If  $\psi$ is a continuous function on
$S_\phi(M)$ such that $\psi$ can be expressed  as a pointwise
limit of a sequence of algebraic combinations of functions of the
form $p\mapsto \phi(p,b)$, $b\in M$, then $\psi$ is called a {\em
$\phi$-definable relation over $M$}. A definable relation
$\psi(y)$ over $M$ defines $p \in S_\phi(M)$ if $\phi(p,b)
=\psi(b)$ for all $b\in M$.

The next notion is more natural and less technically involved
than measure theoretic notion, Definition~\ref{stab1} below. (See
Definition 7.1 in \cite{Ben-Gro}.)

\begin{dfn}
\label{stab2} {\em A formula $\phi(x,y)$ is called {\em stable in
a structure} $\M$ if there are no $r>s$ and infinite sequences
$a_n, b_n \in M$ such that for all $i>j$: $\phi(a_i,b_j)\geq r$
and $\phi(a_j,b_i)\leq s$. A formula $\phi$ is {\em stable in a
theory} $T$ if it is stable in every model of $T$.}
\end{dfn}

It is easy to verify that $\phi(x,y)$ is stable in $\M$ if
whenever $a_n, b_n \in M$ form two sequences, then
$$\lim_n \lim_m \phi(a_n,b_m) = \lim_m \lim_n \phi(a_n,b_m),$$
  provided both limits exist.

\begin{fct}[Grothendieck's Criterion, \cite{Gro}] \label{Fact2}
Let $X$ be an arbitrary topological space, $X_0\subseteq X$ a
dense subset. Then the following are equivalent for a subset
$A\subseteq C_b(X)$:

\begin{itemize}
             \item [{\em (i)}]  The set A is relatively weakly compact in $C_b(X)$.
             \item [{\em (ii)}] The set
$A$ is bounded, and whenever $f_n\in A$ and $x_n\in X_0$ form two
sequences we have $$\lim_n \lim_m f_n(x_m) =\lim_m \lim_n
f_n(x_m),$$  whenever both limits exist.
\end{itemize}
\end{fct}

\subsection{Fundamental theorem of stability}
In \cite{BU}, Ben Yaacov and Usvyatsov proved a continuous
version of the definability of types in a stable theory, which is
a generalization of the classical one.
Roughly speaking, in continuous logic, for a stable formula
$\phi$, the number of $\phi$-types is controlled by the number of
continuous functions on the   space of $\phi$-types. A similar
result holds for a stable formula in integral logic. Also,
another result shows that for an almost dependent formula $\phi$
(see Definition~\ref{stab1} below), the number of $\phi$-types (up
to an equivalence relation) is controlled by the number of
measurable functions  on the space of $\phi$-types.

On the other hand, in \cite{Ben1} and \cite{Ben2}, Ben Yaacov
studied probability algebras and $L^1$-random variables in the
frameworks of  compact abstract theories (cats) and of continuous
logic.  Note that in this paper we shall \textit{not} identify
measurable functions with their class in $L^1$. Thus, in contrast
to \cite{Ben1} and \cite{Ben2}, the theory of a probability
structure is not necessarily stable.

Now, we come quickly to the following theorem. The proof  is
essentially similar to that in \cite{Ben-Gro}, but it works for
measure structures.

\begin{thm}[Definability of types] \label{type-dfn}  Let $\phi(x,y)$ be a formula stable
  in a structure $\M$. Then every $p\in S_\phi(M)$ is definable by
a unique  $\tilde{\phi}$-definable relation $\psi(y)$ over $M$,
where $\tilde{\phi}(y,x)=\phi(x,y)$.
\end{thm}
 {\bf Proof.}  Let $X=S_{{\phi}}(M)$  and let $X_0 \subseteq
X$ be the collection of those types realized in $M$, which is
dense in $X$. Since $X$ is compact,
  the weak topology on $C(X)$
coincides with the topology of pointwise convergence. Since every
formula is bounded, the set $A = \{\phi^a:p\mapsto \phi(a,p)~ |~
a\in M\}\subseteq C(X)$ is bounded.
 By Fact \ref{Fact2}, since $\phi$ is  stable in $M$,
so $A$ is relatively poinwise compact in $C(X)$. Let $p(x)\in
S_\phi(M)$, and let $a_i\in M$ be any net such that
$\lim_i\text{tp}_\phi(a_i/M)=p$. Since $A$ is relatively pointwise
compact, there is a $\psi\in C(X)$ such that
$\lim_i\phi^{a_i}(y)=\psi(y)$.
 By Theorem~8.20 in \cite{KN}, $\psi$ is
the closure point of a sequence $\phi^{a_n}(y)$ of the family
$\{\phi^{a_i}(y)\}_i$, and there is a subsequence
$\phi^{a_{n_k}}(y)$ such that $\lim_k\phi^{a_{n_k}}(y)=\psi(y)$.
Clearly, $\psi(y)$ is a $\tilde{\phi}$-definable relation over
$M$, and for $b\in M$ we have $\phi(p, b)=\lim_k \phi(a_{n_k},b)
= \psi(b)$. Therefore, $p$ is definable by a
$\tilde{\phi}$-definable relation $\psi$ over $M$. If $p$ is
definable by $\psi_1,\psi_2$, then $\psi_1(b)=\psi_2(b)$ for all
$b\in M$. Since $X_0\subseteq X$ is dense, $\psi_1=\psi_2$.~\qed

\bigskip
We are now ready to prove the main theorem of this section.

\begin{cor}[Fundamental Theorem of Stability] \label{Fundamental-Thm}  Let $\phi(x, y)$  be a
formula and $T$ a theory. Then the following are equivalent.
\begin{itemize}
             \item [{\em (i)}] The formula $\phi$ is stable in T.
             \item [{\em (ii)}] For every model ${\M}\vDash T$, every $\phi$-type over $M$ is definable
                                by a $\tilde{\phi}$-predicate over $M$.
             \item [{\em (iii)}] For each cardinal $\lambda=\kappa^{\aleph_0}\geq
                               |T|$, and model ${\M}\vDash T$ with $|M|\leq\lambda$, $|S_\phi(M)|\leq\lambda$.
             \item [{\em (iv)}] There exists a
cardinal $\lambda=\kappa^{\aleph_0}\geq |T|$ such that for every
model ${\M}\vDash T$, $|M|\leq\lambda$ then
$|S_\phi(M)|\leq\lambda$.
\end{itemize}
\end{cor}
{\bf Proof.} We proved (i) implies (ii) in Theorem \ref{type-dfn}.
The implications (ii) $\Longrightarrow$ (iii) $\Longrightarrow$
(iv) are clear. For (iv) $\Longrightarrow$ (i), use many type
argument and the downward L\"{o}wenheim-Skolem theorem
(Proposition~5.13 in \cite{BP}).~\qed

\subsection{Cantor-Bendixson rank}
 Let $\M$ be a structure. By Remark \ref{topomeasure}, ${\N}=(S({\M}),\mu)$ is
an elementary extension of $\M$, and a very unlikely one from the
point of view of classical model theory. Moreover, $\N$ is a
topological measure space, $N$ is compact and $\mu$ is a Radon
measure. Similarly, for a formula $\phi(x,y)$, the structure
${\N}_\phi=(S_\phi({\M}),\mu_\phi)$ also is. In fact, ${\N}_\phi$
has further structures:

\begin{dfn}[\cite{BU}]
 {\em A (compact) \emph{topometric} space is a triplet $\langle X, \tau,d\rangle$,
  where $\tau$ is a
  (compact) Hausdorff topology and $d$ a metric on $X$, satisfying:
  (i) The metric topology refines the topology.
  (ii) For every closed $F \subseteq X$ and $\epsilon > 0$, the closed
    $\epsilon$-neighbourhood of $F$ is closed in $X$ as well.}
\end{dfn}

\begin{fct}
${\N}_\phi$ is a compact topometric space.
\end{fct}
{\bf Proof.} For $p,q\in S_\phi(M)$, define
$d(p,q)=\sup\{|\phi(p,a)-\phi(q,a)|:a\in M\}$. Clearly, $d$ is a
metric on $S_\phi(M)$, and  the topology generated by $d$
sometimes called the {\em uniform topology}. On the other hand, we
know that $p_\alpha\to p$ in the logic topology $\tau$ iff
$\phi^{p_\alpha}\to \phi^p$ in the topology of pointwise
convergence, or equivalently, iff $\phi^{p_\alpha}\to \phi^p$ in
the weak topology. Now, it is easy to verify that
$(S_\phi(M),\tau,d)$ is a compact topometric space.~\qed

\begin{rmk} \label{topometric}
{\em Let $U$ be an Archimedean Riesz space with order unit $e$.
Then it can be embedded as an order-dense and norm-dense Riesz
subspace of $C(X)$, where $X$ is a compact Hausdorff space (see
Fact~\ref{embed to C(X)}). For $a,b\in X$, define
$d(a,b)=\sup\{|f(a)-f(b)|:f\in C(X)\}$. Clearly, $(X,d)$ is a
compact topometric space.
 Therefore, all results in this paper
can be extended to Archimedean Riesz space with order unit, and
our approach is appropriate for continuous logic as well as
operator logics (cf.~\cite{Mofidi}).}
\end{rmk}

We have the following continuous version of the Cantor-Bendixson
rank.

\begin{dfn}[\cite{BU}]  {\em  Let $X$ be a compact topometric space. For a
fixed $\epsilon > 0$, we define a decreasing sequence of closed
subsets $X_{\epsilon,\alpha}$ by induction:
 \begin{align*}
X_{\epsilon,0} &= X \\
 X_{\epsilon,\alpha} &=  \bigcap_{\beta<\alpha}X_{\epsilon,\beta}~~~~~~~~~\mbox{ for $\alpha$ a limit ordinal} \\
X_{\epsilon,\alpha+1} &= \bigcap\{F\subseteq
X_{\epsilon,\alpha}:\mbox{  $F$ is closed and
$\mbox{diam}(X_{\epsilon,\alpha}\setminus F)\leq\epsilon$}\}
\\ X_{\epsilon,\infty}&= \bigcap_\alpha X_{\epsilon,\alpha}.
\end{align*}
Where the {\em diameter} of a subset $U\subseteq X$ is defined
    $$\mbox{diam}(U) = \sup\{d(x,y)\colon x,y \in U\}.$$

For any non-empty subset $U\subseteq X$ we define its {\em
$\epsilon$-Cantor-Bendixson rank in $X$} as:
$$CB_{X,\epsilon}(U)=\sup\{\alpha:U\cap X_{\epsilon,\alpha}\neq \emptyset\}\subseteq
Ord\cup\{\infty\}$$}
\end{dfn}

The next result characterizes stability in terms of CB ranks. We
remark that a structure $\M$ is $\omega$-saturated if every
$1$-type over a finite tuple in $M$ is realized in $\M$.

\begin{pro}[cf.~\cite{BU}] $\phi$ is stable iff for any $\omega$-saturated model ${\M}\vDash T$ where $|M|=(|T|+\kappa)^{\aleph_0}$
we have $CB_{S_\phi(M),\epsilon}(S_\phi(M))<\infty$ for all
$\epsilon$.
\end{pro}
{\bf Proof.} Let $\kappa>|T|$ be any cardinal such that
$\kappa=\kappa^{\aleph_0}$.
 Let $\lambda$ be the least
cardinal such that $2^\lambda>\kappa$. Assume that
$Y=S_\phi(M)_{\epsilon,\infty}$ is nonempty. Therefore $Y$ is
compact and if $U \subseteq Y$ is relatively open and
  non-empty then $\textrm{diam}(U) > \epsilon$.
We can therefore find non-empty open sets $U_0,U_1$ such that
$\bar U_0,\bar U_1 \subseteq U$ and $d(U_0,U_1) > \epsilon$.
  Now, if $p\in U_0,q\in U_1$ then $d(p,q)>\epsilon$.
Proceed by induction. If $\M$ be $2^{<\lambda}$-saturated and
$(2^{<\lambda})^{\aleph_0}=2^{<\lambda}$, then we can find a model
${\M}_0\preceq {\M}$ of cardinality $2^{<\lambda}$ and the types
$\{p_\alpha\}_{\alpha<2^\lambda}\subseteq S_\phi(M_0)$ such that
$d(p_\alpha,p_{\alpha'})>\epsilon$ for all $\alpha\neq \alpha'$.
  Therefore, $\|S_\phi(M_0)\|>|M_0|$, i.e., the density
character of $S_\phi(M_0)$ is bigger than the cardinality of
$M_0$.

  The converse is also standard.~\qed

\subsection{Stability and amenability} Now we return to
analytic concepts.  A topological group is called {\em
precompact} if it is isomorphic to a subgroup of a compact group.
Assume that $G$ acts on a set $X$. A bounded function $f$ on $X$
is called {\em weakly almost periodic}  if the $G$-orbit of $f$
is weakly relatively compact in the Banach space $l^\infty(X)$
of   all bounded real-valued functions on $X$ equipped with the
supremum norm. For a topological group $G$, denote by $WAP(G)$
the space of all continuous weakly almost periodic functions on
$G$.

\begin{fct}
Assume that $G$ is a topological group and its theory, $T_G$,
 is satisfiable. Then the following are equivalent:
\begin{itemize}
             \item [{\em (i)}]  $T_G$ is stable.
             \item [{\em (ii)}] $G$ is precompact.
\end{itemize}
\end{fct}
 {\bf Proof.} We know that $T_G$ is
stable (i.e., $LUC(G)$ is weakly compact) if and only if
$LUC(G)=WAP(G)$. By Theorem 4.5 in \cite{MPU}, $LUC(G)=WAP(G)$ if
and only if  $G$ is precompact.~\qed

\begin{cor}
Assume that $G$  and $T_G$ are as above. If $T_G$ is stable, then
$G$ is uniquely amenable.
\end{cor}
 {\bf Proof.} It is known that for every precompact group $G$, the algebras $LUC(G)$
and $LUC(\widehat G)$ are canonically isomorphic, where $\widehat
G$ denotes the compact completion of $G$. Also, every compact
group provides an obvious example of a uniquely amenable group,
for which the unique invariant mean comes from the Haar measure.
So $G$ is uniquely amenable since $\widehat G$ is.~\qed

\section{NIP} \label{NIP}
 In \cite{Talagrand}, Talagrand gave the first explicit  definition
of stable set of functions. In fact, the notion of stable set of
functions (\cite[465B]{Fremlin4}) is a measure-theoretic version
of a well-known model-theoretic property, the dependence property.
The definition is not obvious, but given this the basic
properties of stable sets listed in \cite[465C]{Fremlin4} are
natural and easy to check, and we come quickly to the fact that
(for complete locally determined spaces) pointwise bounded stable
sets are relatively pointwise compact sets of measurable functions
(Fact~\ref{Fact1}).  We are now ready for the main definition
which is an adapted version of Definition 465B in \cite{Fremlin4}.

\subsection{Almost dependence property}
\begin{dfn} \label{stab1}
 {\em A formula $\phi(x,y)$
has the  {\em almost dependence property}, or is {\em almost
dependent}, in a structure $\M$ if the set
$A=\{\phi(x,b),\phi(a,y):a,b\in M\}$ is a stable set of functions
in the sense of Definition~465B in \cite{Fremlin4}, that is,
whenever $E\subseteq M$ is measurable, $\mu(E)>0$ and $s<r$ in
$\mathbb{R}$, there is some $k\geq 1$ such that
$(\mu^{2k})^*D_k(A, E,s,r)<(\mu E)^{2k}$ where
\begin{align*}
D_k(A, E,s,r) = \bigcup_{f\in A}\{w\in & E^{2k}:f(w_{2i})\leq s,
~f(w_{2i+1})\geq r  \textrm{ for } i<k\}.
\end{align*}
 A formula $\phi$ has the {\em almost dependence property}
 in a theory $T$ if it has the almost dependence property in
every model of $T$.}
\end{dfn}

\begin{note} {\em Assume that for each $s< r$ and $k\in
\mathbb{N}$ the set $D_k(A, E,s,r)$ is measurable in $\M$. Then
it is easy to verify that  $\phi(x,y)$ fails to be almost
dependent in $\M$ if and only if there exist $E\subseteq M$, with
$\mu(E)>0$ and $s<r$ in $\mathbb{R}$, such that for each $k\geq
1$, and almost each $w\in E^k$, for each
$I\subseteq\{1,\ldots,k\}$, there is $f\in A$ with $f(w_i)\leq s$
for $i\in I$ and $f(w_i)\geq r$ for $i\notin I$ (see
\cite{Talagrand}, Proposition~4).

In the above definition if $\mu(E)\geq\epsilon>0$ then we say that
$\phi$ fails to be almost $\epsilon$-dependent, or it has the
{\em $\epsilon$-FD property}. It is an easy exercise to show that
the $\epsilon$-FD property is a {\em first order} property (in
integral logic), or equivalently it is expressible. Clearly,
$\phi$ has the almost dependence property if it has not the
$\epsilon$-FD property for all $\epsilon>0$.

 Note that the sets $A_1=\{\phi(a,y):a\in M\}$ and $A_2=\{\phi(x,b):b\in M\}$ are
dependent if and only if $A=A_1\cup A_2$ is dependent  (cf.
Proposition~465C(a),(d) in \cite{Fremlin4}).
 On the other hand, one can easily define the (exact)
{\em dependence property}. For this, we say $\phi$ fails to be
{\em dependent}, or is {\em independent}, in $\M$ iff  there exist
$s<r$ in $\mathbb{R}$, such that for each $k\geq 1$, there are
$w_1,\ldots,w_k\in M$, such that for each
$I\subseteq\{1,\ldots,k\}$, there is $f\in A$ with $f(w_i)\leq s$
for $i\in I$ and $f(w_i)\geq r$ for $i\notin I$.
 Clearly, a dependent
formula (or theory) is necessarily almost dependent.}
\end{note}

We come quickly to the following fact which is an adapted version
of Proposition 465D in \cite{Fremlin4}.

\begin{fct} \label{Fact1}
Let ${\M}=(M,\Sigma,\mu)$ be a structure such that $\mu$ is a
complete locally determined measure space, and $\phi(x,y)$ an
almost dependent formula.  Since every formula is bounded, so
$\phi$ is. Therefore, $A=\{\phi(x,b),\phi(a,y) : a,b\in M\}$ is
relatively compact in the space of measurable functions for the
topology of pointwise convergence.
\end{fct}

First we compare our notions.

\begin{pro}
Let $\phi(x,y)$  be a stable formula in a theory $T$. Then $\phi$
is almost dependent  in $T$.
\end{pro}
{\bf Proof.}   Assume that $\phi$ fails to be almost dependent.
Therefore, there is a model ${\M}\vDash T$, $E\subseteq M$, with
$\mu(E)>0$, and $r>s$ in $\mathbb R$ such that
$(\mu^{2k})^*D_k(A, E,s,r)=(\mu E)^{2k}$ for each $k$. Then it is
easy to verify that for each $k$ there are  finite sequences
$a_n, b_n \in E$, $n\leq k$ such that for all $j<i\leq k$:
$\phi(a_i,b_j)\geq r$ and $\phi(a_j,b_i)\leq s$.
 Now, by the
compactness theorem of model theory, there is an elementary
extension ${\N}\succ{\M}$ such that $\phi$ is not stable in $\N$.
Thus, $\phi$ is not stable in $T$.~\qed

\bigskip
 To summarize:
$$
  \phi \textrm{ is stable}  \ \ \   \Longrightarrow \ \ \   \phi \textrm{ is dependent} \ \ \  \Longrightarrow \ \ \  \phi \textrm{ is almost dependent} \\
$$

By a result of Bourgain,  Fremlin and  Talagrand
\cite[Theorem~2F]{BFT}, one can easily check that a formula $\phi$
is (exactly) dependent if and only if it is almost dependent for
each Radon measure. We will study their connection in a future
work.

\subsection{Almost definability of types}
Here, a result similar to stable case can be proved for the almost
dependence property. For this, we need some definitions. Let
$\psi$ be a measurable function on $(S_\phi(M),\mu_\phi)$ where
$\mu_\phi$ is the unique Radon measure induced by
$\phi^{\M}(x,b)$ for all $b\in M$. Then $\psi$ is called an {\em
almost $\phi$-definable relation over $M$} if there is a sequence
$g_n:S_\phi(M)\to \mathbb{R}$, $|g_n|\leq |\phi|$, of continuous
functions such that $\lim_n g_n(b)=\psi(b)$ for almost all $b\in
S_\phi(M)$. (We note that by the Stone-Weierstrass theorem every
continuous function $g_n:S_\phi(M)\to \mathbb{R}$ can be
expressed as a uniform limit of a sequence of algebraic
combinations of functions of the form $p\mapsto \phi(p, b)$,
$b\in M$.)  An almost definable relation $\psi(y)$ over $M$
defines $p \in S_\phi(M)$ if $\phi(p,b) =\psi(b)$ for almost all
$b\in M$, and in this case we say that $p$ is {\em almost
definable}. Assume that every type $p$ in $S_\phi(M)$ is almost
definable by a measurable function $\psi^p$.
 Then, we say that $p$ is {\em almost equal to $q$},
denoted by $p\equiv q$, if $\psi^p(b)=\psi^q(b)$ for almost all
$b\in M$. Define
 $[p]=\{q\in S_\phi(M): p\equiv q\}$ and  $[S_\phi](M)=\{[p]:p\in S_\phi(M)\}$.

\begin{thm}[Almost definability of types]  \label{almost-Definable-Thm}  Let $\phi(x,y)$ be a formula almost dependent
  in a structure $\M$. Then every $p\in S_\phi(M)$ is almost definable by
a (unique up to measure) almost $\tilde{\phi}$-definable relation
$\psi(y)$ over $M$, where $\tilde{\phi}(y,x)=\phi(x,y)$.
\end{thm}
 {\bf Proof.}
We know that $({\M},\mu^{\M}_\phi)\preceq(S_\phi(M),\mu_\phi)$.
First, we assume that $(S_\phi(M),\mu_\phi)$ is minimal, i.e.,
$\mu_\phi$ is Baire and it is not nesessarily  Radon. (One can
easily verify that the subspace measure
${\mu_\phi}\upharpoonright_M$ is the measure $\mu^{\M}_\phi$.
Therefore,
 by Proposition 465C(n) in \cite{Fremlin4}, since the set $\{\phi(a,y)\upharpoonright_M: a\in M\}\subseteq \mathbb{R}^{M}$
  is almost dependent with respect to  the subspace measure ${\mu_\phi}\upharpoonright_M$,
   the set $A=\{\phi(a,y): a\in M\}\subseteq \mathbb{R}^{S_\phi(M)}$ is also almost dependent with respect to  $\mu_\phi$.)
   By Proposition 465C(i) in \cite{Fremlin4}, the set $A$ is also
almost dependent with respect to the completion  $\hat\mu_\phi$ of
$\mu_\phi$. Now, let $p(x)\in S_\phi(M)$, and let $a_i\in M$ be
any net such that $\lim_i\text{tp}_\phi(a_i/M)=p$. Since
$\hat\mu_\phi$ is complete, by Fact \ref{Fact1},
 there is a $\hat\mu_\phi$-measurable function $\psi$
such that $\lim_i\phi^{a_i}(y)=\psi(y)$. Let $\bar\mu_\phi$ be
the unique   extension of $\mu_\phi$ to a Radon measure. Thereby
it is also an extension of $\hat\mu_\phi$. Since $\bar\mu_\phi$
is Radon, by Proposition 7.9 in \cite{Folland},
 there is a sequence $g_n$ of continuous functions  on $S_\phi(M)$
such that $g_n\to \psi$ in $L^1(\bar\mu_\phi)$, and hence by
Corollary 2.32 in \cite{Folland} a subsequence (still denoted by
$g_n$) that converges to $\psi$ $\bar\mu_\phi$-a.e. Clearly,
$\psi$ is unique up to the measure   $\bar\mu_\phi$.~\qed

\begin{cor}  \label{almost-Fundamental-Thm}  Let $\phi(x, y)$  be a
formula and $T$ a theory. Then (i) $\Longrightarrow$ (ii)
$\Longrightarrow$ (iii).
\begin{itemize}
             \item [{\em (i)}] The formula $\phi$ is almost dependent in T.
             \item [{\em (ii)}] For every model ${\M}\vDash T$, every $\phi$-type over $M$ is almost definable
                                by a $\tilde{\phi}$-predicate over $M$.
             \item [{\em (iii)}] For each cardinal $\lambda=\kappa^{\aleph_0}\geq
                               |T|$, and model ${\M}\vDash T$ with $|M|\leq\lambda$, $|[S_\phi](M)|\leq\lambda$.
\end{itemize}
\end{cor}
{\bf Proof.}  Clear.~\qed

\subsection{Almost Cantor-Bendixson rank}
A  result similar to the Cantor-Bendixson rank for stable formulas
holds for the almost dependence property. For this we need some
definitions. For a $\mu_\phi$-measurable function
$\xi:S_\phi(M)\to[-\flat_\phi,\flat_\phi]$ where $\flat_\phi$ is
the universal bound of $\phi$, let
$$[\xi]=\big\{\chi:S_\phi(M)\to [-\flat_\phi,\flat_\phi]~\big|~ \chi
\mbox{ is $\mu_\phi$-measurable and } \chi=\xi~ a.e\big\}.$$
 Let $L^1_{\phi}=\{[\xi]~|~\xi:S_\phi(M)\to[-\flat_\phi,\flat_\phi] \mbox{ is $\mu_\phi$-measurable}\}$.
We show that $L^1_\phi$ has a natural  compact topometric
structure. Indeed, let ${\frak
d}([\xi],[\xi'])=\int|\xi-\xi'|d\mu_\phi$, and
 $[\xi_\alpha]\to^{\frak T} [\xi]$ iff $I([\xi_\alpha])\to I([\xi])$ for all
 $I\in (L^1)^*$. In fact, the topology generated by the metric ${\frak d}$ is the norm topology on $L^1$ and ${\frak T}$ is
 the weak topology generated by $(L^1)^*$.
 Now, it is easy to verify that $(L^1_{\phi},{\frak d},{\frak T})$ is a compact topometric
 space. Indeed, since $L_\phi^1$ is uniformly integrable, by
 Theorem 247C in \cite{Fremlin2}, $L_\phi^1$ is relatively weakly
 compact. Also, $L_\phi^1$ is closed in the norm topology. It is
 well-known that for a convex subset of a locally convex space, the weak
 closure is equal to the norm closure. Therefore, $L_\phi^1$ is weakly
 closed, and hence it is weakly compact.
 On the other hand, it is well-known that the norm $L^1$ is weakly lower
 semicontinuous (cf. Lemma 6.22 in \cite{Ali}). To summarize, $L_\phi^1$ is a compact topometric space.

We remark that if the types $p,q$ are definable by measurable
functions $\psi^p,\psi^q$, then $p\equiv q$ iff
$\psi^p(b)=\psi^q(b)$ for almost all $b\in M$, or equivalently,
iff $[\psi^p]=[\psi^q]$. (Note  that since ${\M}\preceq
(S_\phi(M),\mu_\phi)$ therefore $\psi^p(b)=\psi^q(b)$ for almost
all $b\in M$ iff $\psi^p=\psi^q$ $\mu_\phi$-almost everywhere.)
  Therefore, if $|M|=\kappa^{\aleph_0}$ and $\phi$ is a
formula almost dependent in the structure $\M$, then
$|L^1_{\phi}|=|[S_\phi](M)|=\kappa^{\aleph_0}$. (See Theorem
\ref{almost-Definable-Thm} and the definitions before it.)
Thereby:

\begin{pro} If $\phi$ is almost dependent then for any   $\omega$-saturated model ${\M}\vDash T$ where $|M|=(|T|+\kappa)^{\aleph_0}$
we have $CB_{L^1_{\phi},\epsilon}(L^1_{\phi})<\infty$ for all
$\epsilon$.
\end{pro}

Almost dependence property is linked with a  notion  of another
area. Historically, this property arose naturally when Talagrand
and Fremlin were studying pointwise compact sets of measurable
functions; they found that in many cases a set of functions was
relatively pointwise compact because it was almost dependent.
Later did it appear that the concept was connected with {\em
Glivenko-Cantelli classes} in the theory of empirical measures,
as explained in \cite{Talagrand}. Also, a version of {\em
Vapnik-Chervonenkis dimension} which is suitable for measure
structures can be defined, and will be studied in a future work.

\section{Conclusion}
In the first part of this paper we studied some concrete analytic
structures. This study led us to the natural and correct notion
of types. The perspective on types in this paper can be used in
other logics. For example, this approach seems to be appropriate
for continuous logic \cite{BU} as well as operator logics
\cite{Mofidi}. Note that by Remark \ref{topometric},
 every Archimedean Riesz space with order unit admits a natural
 compact topometric structure. Therefore, the most of results in this paper can be extended to
Archimedean Riesz spaces. Also, the notion of forking and
independence, and their connections to measure theory can be
studied. On the other hand, one can do much more classifications,
 the strict order property and others. We will study them
elsewhere. Finally, all these results suggest that many
interesting analytic concepts may be studied by model theoretic
methods. Also, these methods provide a new view on the related
subjects in Analysis, and open some fruitful areas of research on
similar questions.

\bigskip\noindent
 {\bf Acknowledgements.}
I want to thank John T. Baldwin and Massoud Pourmahdian for their
interest in reading a copy of this article and for their comments
and corrections. I thank the anonymous referee for his detailed
suggestions and corrections; they helped to improve significantly
the exposition of this paper.  I would like to thank the Institute
for Basic Sciences (IPM), Tehran, Iran.  Research partially
supported by IPM grants no. 91030032 and 92030032.



\begin{thebibliography}{EGGN07} \label{ref}
\bibitem[AB06]{Ali} C. D. Aliprantis and K. C. Border, {\em Infinite Dimensional Analysis, A Hitchhiker's Guide}  (3rd ed.), Springer-Verlag, Berlin Heidelberg (2006).
\bibitem[BP09]{BP} S. M. Bagheri and M. Pourmahdian,  The logic of integration, {\em Arch. Math. Logic.} {\bf 48},  465-492 (2009).
\bibitem[Ben14]{Ben-Gro}  I. Ben-Yaacov,  Model theoretic stability and definability of types. after A. Grothendiek, Bull. Symbolic Logic, 20 (2014), 491-496.
\bibitem[Ben13]{Ben2} I. Ben-Yaacov, On theories of random variables, Israel J. Math. 194 (2013), no. 2, 957-1012
\bibitem[Ben06]{Ben1} I. Ben-Yaacov,   Schr\"{o}dinger's cat, Israel Journal of Mathematics \textbf{153} (2006), 157–191.
\bibitem[BBHU08]{BBHU} I. Ben-Yaacov, A. Berenstein, C. W. Henson, A. Usvyatsov, \emph{Model theory for metric structures},
     Model theory with Applications to Algebra and Analysis, vol. 2
   (Z. Chatzidakis, D. Macpherson, A. Pillay, and A. Wilkie, eds.),
    London Math Society Lecture Note Series, vol. 350, Cambridge University Press, 2008.
\bibitem[BU10]{BU} I. Ben-Yaacov, A. Usvyatsov, Continuous first order logic and local stability, Transactions of the American Mathematical Society \textbf{362} (2010), no. 10, 5213-5259.
\bibitem[BFT78]{BFT} J. Bourgain, D. H. Fremlin, and M. Talagrand. Pointwise compact sets of baire-measurable functions. American Journal of Mathematics, 100(4):pp. 845–886, 1978.
 \bibitem[CK66]{CK} C. C. Chang,  H. J. Keisler, {\em Continuous Model Theory}. Princeton University Press, Princeton (1966).
\bibitem[DS58]{Dunford} N. Dunford and J. T. Schwartz, {\em Linear operators}, Part I, Pure and Appl. Math. No. 7, Interscience, New York, 1958.
\bibitem[Fol99]{Folland} G. B. Folland, {\em Real Analysis, Modern Techniquess and Their Appilcations}, 2nd ed. Wiley, New York (1999).
\bibitem[Fre03]{Fremlin2} D. H. Fremlin, {\em Measure Theory}, vol.2, (Broad Foundations, Torres Fremlin, Colchester, 2003).
\bibitem[Fre04]{Fremlin3} D. H. Fremlin, {\em Measure Theory}, vol.3, (Measure Algebras, Torres Fremlin, Colchester, 2004).
\bibitem[Fre06]{Fremlin4} D. H. Fremlin, {\em Measure Theory}, vol.4, (Topological Measure Spaces, Torres Fremlin, Colchester, 2006).
\bibitem[Gro52]{Gro} A. Grothendieck, Crit\`{e}res de compacit\'{e} dans les espaces fonctionnels g\'{e}n\'{e}raux, American Journal of Mathematics {\bf 74} ,168-186, (1952).
\bibitem[HC75]{Herer} W. Herer and J.P.R. Christensen,  On the existence of pathological submeasures and the construction of exotic topological groups, {\em Math. Ann}. {\bf 213},  203–210, (1975).
\bibitem[HR63]{Hewitt} E. Hewitt and K. A. Ross, Abstract Harmonic Analysis, Vol. I, Springer-Verlag, Berlin (1963).
\bibitem[Hoo78]{Hoover} D. N. Hoover, Probability logic, {\em Ann.  Math. Log.} {\bf 14}, 287-313 (1978).
\bibitem[Kei85]{Keisler} H. J. Keisler, Probability quantifiers, In: J. Barwise, S. Feferman, (eds.) {\em Model Theoretic Logic}, Springer, Berlin (1985).
\bibitem[KN63]{KN}  J. L. Kelley, I. Namioka, {\em Linear topological spaces},  New York (1963).
\bibitem[KA12]{KA} K. Khanaki and  M. Amini, Haar measure and integral logic, {\em  Math. Log. Quart.} {\bf 58}, No. 4-5, 294-302 (2012).
\bibitem[KB11]{KB} K. Khanaki and S. M. Bagheri, Random variables and integral logic, {\em  Math. Log. Quart.} {\bf 57}, No. 5, 494-503 (2011).
\bibitem[MPU01]{MPU} M. Megrelishvili, V. G. Pestov,  V. V. Uspenskij,  A note on the precompactness of weakly almost periodic
groups, In: Nuclear groups and Lie groups, E. Martin Peinador, J.
Nunez Garcia (eds.), Heldermann Verlag, 2001 (Research and
Exposition in Mathematics, Volume 24), pp. 209-216.
\bibitem[Mit70]{Mitch2} T. Mitchell, Topological semigroups and fixed points, {\em Illinois J. Math}., {\bf 14} (1970), 630-641.
\bibitem[Mof12]{Mofidi} A. Mofidi, Quantified universes and ultraproducts, {\em  Math. Log. Quart.} 58.1-2 (2012): 63-74.
\bibitem[Nam67]{Namioka} I. Namioka,  On certain actions of semi-groups on $L$-spaces, {\em Studia Math}., vol. 29, pp. 63-77, (1967).
\bibitem[Pat86]{Paterson} A. L. T. Paterson,  Nonamenability and Borel paradoxical decompositions for locally compact groups. {\em Proc. Amer. Math. Soc}.  {\bf 96}, 89-90, (1986).
\bibitem[Pes99]{Pestov1} V. Pestov,  Topological groups: where to from here?, {\em  Topology Proc}. {\bf 24},  421-502, (1999).
\bibitem[Pes02]{Pestov2} V. Pestov,  Ramsey-Milman phenomenon, Urysohn metric spaces, and extremely amenable groups,  {\em Israel J. Math}. {\bf 127}, 317–357, (2002).
\bibitem[Roy88]{Royden} H. L. Royden, {\em Real Analysis} (3rd ed.), Macmillan, New York (1988).
\bibitem[Run02]{Runde} V. Runde, {\em Lectures on Amenability}, Lecture Notes in Math. 1774, Springer-Verlag Berlin Heidelberg (2002).
\bibitem[Sim14]{S} P. Simon, Rosenthal compacta and NIP formulas, arXiv 2014
\bibitem[Tal87]{Talagrand} M. Talagrand, The Glivenko-Cantelli problem, {\em Ann. of Probability} 15 (1987) 837-870.
\bibitem[Wag85]{Wagon} S. Wagon, {\em The Banach-Tarski Paradox}. Cambridge University Press, (1985).
\end{thebibliography}
\end{document}